\documentclass[preprint, 1p, sort&compress, a4paper]{elsarticle}

\journal{Engineering Analysis with Boundary Elements}

\usepackage{amsmath}
\usepackage{amsfonts}
\usepackage{amssymb}
\usepackage{mathrsfs}	
\usepackage{bm}
\usepackage{fancybox}
\usepackage{makeidx}
\usepackage{graphicx}
\usepackage{subfigure}
\usepackage{epsfig}
\usepackage{verbatim}
\usepackage{hyperref}
\usepackage{cases}
\usepackage{booktabs}
\usepackage{upgreek}
\usepackage{tabularx,multirow}
\usepackage{color}
\usepackage[active]{srcltx}
\usepackage{citesort}
\usepackage{ulem}   
\let \oldmarginpar \marginpar
\renewcommand{\marginpar}[1]{\oldmarginpar{\color{red}{#1}}}

\newcounter{algleo}
\newlength{\lefttab}
\newlength{\numberoffset}
  {\trivlist
   \topsep=0pt\parsep=0pt\itemsep=0pt
   \addtolength{\lefttab}{1.25em}
   \addtolength{\numberoffset}{1.25em}
   \leftskip=\lefttab}%
  {\endtrivlist}

\definecolor{MyGray}{rgb}{0.5, 0.5, 0.5}
\definecolor{MyDarkGreen}{rgb}{0.0, 0.46, 0.29}
\definecolor{MyLightBlue}{rgb}{0.85, 0.85, 1.0}
\definecolor{MyLightGreen}{rgb}{0.9, 1.0, 0.9}

\begin{document}

\begin{frontmatter}

\title{A fast directional BEM for large-scale acoustic problems
	based on the Burton-Miller formulation}

\author[sanpu]{Yanchuang Cao}
\ead{caoyanch@126.com}

\author[sanpu,npuicma]{Lihua Wen\corref{cor}}
\ead{lhwen@nwpu.edu.cn}

\author[sanpu,npuicma]{Jinyou Xiao}
\ead{xiaojy@nwpu.edu.cn}

\author[npuicma,uc]{Yijun Liu}
\ead{Yijun.Liu@uc.edu}

\cortext[cor]{Corresponding author}
\address[sanpu]{School of Astronautics, Northwestern Polytechnical 
	University, Xi'an 710072, P. R. China}
\address[npuicma]{Institute for Computational Mechanics and Its 
	Applications, Northwestern Polytechnical University, Xi'an 710072, 
	P. R. China}
\address[uc]{Department of Mechanical Engineering, University of 
	Cincinnati, Cincinnati, Ohio 45221-0072, U.S.A.}

\begin{abstract}
In this paper, a highly efficient fast boundary element method (BEM) for
solving large-scale engineering acoustic problems in a broad frequency 
range is developed and implemented. The acoustic problems are modeled by 
the Burton-Miller boundary integral equation (BIE), thus the fictitious 
frequency issue is completely avoided. The BIE is discretized by using 
the Nystr\"om method based on the curved quadratic elements, leading to 
simple numerical implementation (no edge or corner problems) and high 
accuracy in the BEM analysis. The linear systems are solved iteratively 
and accelerated by using a newly developed kernel-independent wideband 
fast directional algorithm (FDA) for fast summation of oscillatory 
kernels. In addition, the computational efficiency of the FDA is further 
promoted by exploiting the low-rank features of the translation matrices, 
resulting in two- to three-fold reduction in the computational time of 
the multipole-to-local translations. The high accuracy and nearly linear 
computational complexity of the present method are clearly
demonstrated by typical examples. An acoustic scattering problem with
dimensionless wave number $kD$ (where $k$ is the wave number and $D$ is 
the typical length of the obstacle) up to 1000 and the degrees of freedom 
up to 4 million is successfully solved within 10 hours on a computer with 
one core and the memory usage is 24 GB.
\end{abstract}


\end{frontmatter}

\section{Introduction}

The boundary element method (BEM), besides its wide applications in other
branches of science and engineering, has been an important numerical
method in acoustics. This is mainly due to its unique advantages in
dimension reduction, solution accuracy and treating infinite and
semi-infinite domain problems where the radiation condition at infinity
is automatically satisfied. The capability of the acoustic BEM has been
further improved when it is combined with the famous Burton-Miller
formulation\cite{BurtonMiller},
in which the annoying non-uniqueness problem of the conventional Helmholtz
boundary integral equation (BIE) is completely circumvented.

The traditional BEM, however, can not be used for large-scale numerical
simulation, because of the densely populated system matrices and thus the
square scaling of computational cost with respect to $N$, the
degrees-of-freedom (DOF). Fortunately, this limitation has now been
removed to a large extent by various BEM acceleration techniques emerged
in the past three decades. The representative examples are the fast multipole
method (FMM) \cite{fmm, pwfmm, fmm_Darve2000}, $\mathcal{H}$-matrix
\cite{Hmatrix}, wavelet compression \cite{wbem, xwt09}, pre-corrected FFT
\cite{pfft, xy12}, ACA \cite{aca, aca_BR2003}, etc. Of those acceleration
techniques, it is the FMM that has found the most substantial applications
in many areas \cite{n02}, such as acoustics\cite{sl07}, electromagnetics,
elastodynamics \cite{cb08}, to name a few.

The original FMM was first proposed by Greengard and Rokhlin \cite{fmm}
to accelerate the evaluation of interactions of large ensembles of
particles governed by Laplace equation. Its success hinges on the
observation that the interaction via the kernel between well-separated
sets of points is approximately of low rank. In the low frequency regime,
the low rank property still holds for the Helmholtz kernels, thus the FMM
for the Laplace equation is applicable to the Helmholtz equation with
slight modifications \cite{n02}. The complexity of the resulting low
frequency FMM is of order $O(N)$ for a given accuracy. In
the high frequency regime, however, the situation is drastically
different as the low rank assumption is not valid any more.  In fact,
the approximate rank of the translation matrix grows linearly with
the size of the point sets (in terms of the wavelength). However,
Rokhlin \cite{diagonalform} observed that the translation matrix
between well-separated point sets, though of large rank, can be
expressed in ``diagonal form'' in Fourier space. This observation
leads to a high frequency FMM with $O(N\log N)$ complexity.
Nevertheless, this algorithm is only suitable for relatively high
frequency; it becomes numerically unstable in the low frequency regime.
The desirability for a FMM that is accurate and efficient in a broad
frequency range, from low frequency to high frequency, give rise to
the wideband FMM in \cite{wbfmm, wbfmm_Gumerov2009}, which are
essentially hybrid schemes of the FMMs in low and high frequency
regimes. However, besides the success in dealing with both
non-oscillatory and oscillatory kernels, one should note that, the
implementation of the original FMM is highly technical, mathematically
involved and kernel-dependent.  This poses a severe limitation on its
engineering application and its generalization to other problems.

With this background, the kernel-independent FMM (KIFMM) has gain its
popularity in recent years. The main advantage of the KIFMM lies in that
it requires no complicated analytic expansions for the kernels, but only
uses kernel evaluations in its implementation; see e.g. \cite{kifmm,
bbFMM, zh11} and the references therein. For low frequency problems,
we would like to mention two KIFMM algorithms proposed by Ying, Biros
and Zorin in \cite{kifmm} and Fong and Darve in \cite{bbFMM},
both of which have attracted much attention in the past several years.
In \cite{kifmm} the analytic expansions of the original FMM is
replaced with equivalent density representations, while in \cite{bbFMM}
they are replaced with Chebyshev polynomial interpolations.
For applications of the algorithm in \cite{kifmm}, see e.g., 
\cite{bcr11, ybz06}; recent developments of the algorithm in \cite{bbFMM} 
can be found in \cite{opm2l, fda_Chebyshev}.

The KIFMM for high frequency problems are built upon the directional
low rank property of the oscillatory kernels \cite{b91}. Two
representative examples along the this line are the fast directional
algorithm (FDA) proposed by Engquist and Ying \cite{fda}
which uses equivalent density representations (as in \cite{kifmm}),
as well as the one proposed by Messner, Schanz and Darve
\cite{fda_Chebyshev} which uses Chebyshev polynomial interpolations
(as in \cite{bbFMM}). The readers are referred to \cite{fda_Chebyshev}
for the restrictions and advantages of the FDA, and to
\cite{opm2l} for its efficiency enhancement.

This paper will concentrate on the FDA in \cite{fda}. The good accuracy and
efficiency ($O(N \log N)$ complexity) of the algorithm have
been proved theoretically and numerically in \cite{fda}. Furthermore,
the algorithm can be adapted to handle kernels other than the
Helmholtz kernel quite easily, which is not true for most other
existing algorithms. The algorithm was used to accelerate the BEM for
electromagnetic scattering problems by Tsuji and Ying in \cite{ty11}.
However, to the best of the authors' knowledge, no reported work is
devoted to the application of the algorithm to accelerating BEM for
the Burton-Miller BIE, which is essential for solving real-world
acoustic problems.

There are two contributions in this paper. The first one
is to adapt the FDA in \cite{fda} to accelerate the BEM (the resulting
method will be called fast directional BEM, or FDBEM for short,
from now on) for the Burton-Miller BIE; see Section \ref{s:fda}.
The advantages of the FDBEM for Burton-Miller BIE are the following:
(1) it does not require complicated analytical expressions
for the far-field approximations, but only uses kernel evaluations,
which results in a simpler implementation and ease of use; (2) the
algorithm can handle the kernels of the Burton-Miller formulation
in a similar manner.
Therefore, it is more user-friendly than the other FMM algorithms.

The second contribution is achieved by further acceleration for the
translations in the FDA. It is based on the observation that
generally the ranks of most of the M2L translation matrices are much
lower than their dimensions.
The improvement of the algorithm is implemented by the following two steps:
(1) compressing all the translation matrices
into more compact forms, and (2) approximating the compressed M2L matrices
by low rank representations, see Section \ref{s:improve}. Numerical experiments
show that the computational time for M2L translations can be reduced by a
factor of 2 to 3, and
that for the upward and downward passes can be reduced by about 40\%, the memory
consumption is reduced by about 30\%, see Section \ref{s:ne}.

\section{Burton-Miller formulation and its Nystr\"{o}m discretization}
\label{s:nystrom}

\subsection{Acoustic Burton-Miller formulation} \label{BEM:B-M}

The time harmonic acoustic waves in a homogenous and isotropic acoustic medium
$\Omega$ is described by the following Helmholtz equation
\begin{equation}\label{eq:helmholtz}
    \nabla^2 u(\bm{x})+k^2 u(\bm{x})=0,\quad \forall \bm{x} \in \Omega,
\end{equation}
where, $\nabla^2$ is the Laplace operator, $u(\bm{x})$ is the velocity potential
at the point $\bm{x}=(x_1,\,x_2,\,x_3)$ in the physical coordinate system,
$k=\omega/c$ is the wave number, with $\omega$ being the angular frequency
and $c$ being the sound speed. The boundary conditions can be any combinations
of the Dirichlet, Neumann or Robin boundary conditions. For exterior
problems, the Sommerfeld condition at infinity have to be satisfied as well.

By using Green's second theorem, the solution
of Eq. \eqref{eq:helmholtz} can be expressed by integral representation
\begin{equation}\label{eq:IE}
	u(\bm{x}) + \int_\Gamma \frac{\partial G}
	{\partial \bm{n_y}}(\bm{x}, \bm{y}) u(\bm{y}) \mathrm{d} \bm{y} = \int_\Gamma
	G(\bm{x}, \bm{y}) q(\bm{y})
	\mathrm{d} \bm{y} + u^\text{inc}(\bm{x}), \quad \forall \bm{x} \in \Gamma,
\end{equation}
where, $\bm{x}$ denotes the field point and $\bm{y}$ denotes the source point on the boundary
$\Gamma$; $\bm{n_y}$ denotes the unit normal vector at the source
point $\bm{y}$; $q(\bm{y})=\partial u(\bm{y})/\partial \bm{n_y}$ is the
normal gradient of velocity potential, i.e., the normal vibrating velocity
on $\bm{y}$. The incident wave $u^\text{inc}(\bm{x})$ will not be presented
for radiation problems. The three-dimensional (3D) fundamental solution $G$ is given as
\begin{equation}\label{eq:fundsol}
    G(\bm{x},\bm{y})=\frac{e^{\text{i}kr}}{4\pi r},
\end{equation}
with $r = \left| \bm{x}-\bm{y} \right|$ being the Euclidean distance between
the source and the field points, and $\text{i}=\sqrt{-1}$ being the imaginary unit.

Before presenting the BIEs it is convenient to introduce the associated
single, double, adjoint and hypersingular layer operators which are denoted
by $\mathcal{S}$, $\mathcal{D}$, $\mathcal{M}$ and $\mathcal{H}$, respectively;
that is,
\begin{subequations} \label{eq:operator}
\begin{align}
	\mathcal{S} q(\bm{y}) &= \int_\Gamma G(\bm{x}, \bm{y}) q(\bm{y}) \mathrm{d}
		\bm{y}, \label{eq:single}\\
	\mathcal{D} u(\bm{y}) &= \int_\Gamma \frac{\partial G}{\partial \bm{n_y}}
		(\bm{x}, \bm{y}) u(\bm{y}) \mathrm{d} \bm{y}, \label{eq:double}\\
    \mathcal{M} q(\bm{y}) &= \int_\Gamma \frac{\partial G}{\partial \bm{n_x}}
		(\bm{x}, \bm{y}) q(\bm{y}) \mathrm{d} \bm{y}, \label{eq:adjoint}\\
    \mathcal{H} u(\bm{y}) &= \int_\Gamma \frac{\partial^2 G}{\partial \bm{n_x}
		\partial \bm{n_y}} (\bm{x}, \bm{y}) u(\bm{y}) \mathrm{d} \bm{y}. \label{eq:hypersingular}
\end{align}
\end{subequations}
The operator $\mathcal{S}$ is weakly singular and the integral is
well-defined, while the operators $\mathcal{D}$ and $\mathcal{M}$ are
defined in Cauchy principal value sense (CPV). The operator
$\mathcal{H}$, on the other hand, is hypersingular and unbounded
as a map from the space of smooth functions on $\Gamma$ to itself.
It should be interpreted in the Hadamard finite part sense (HFP).
Denoting a vanishing neighbourhood surrounding $\bm{x}$ by
$\Gamma_{\varepsilon}$, the CPV and HFP integrals are those after
extracting free terms from a limiting process to make
$\Gamma_{\varepsilon}$ tends to zero in deriving BIEs \cite{Guiggiani}.

If the field point $\bm{x}$ approaches the boundary $\Gamma$, Eq. \eqref{eq:IE} becomes the conventional BIE (CBIE)
\begin{equation} \label{eq:CBIE}
    c(\bm{x})u(\bm{x})+\mathcal{D} u(\bm{x}) =\mathcal{S} q(\bm{x})+u^\text{inc}(\bm{x}), \quad \bm{x}\in \Gamma,
\end{equation}
where, $c(\bm{x})$ is the free term coefficient which equals to
$1/2$ on smooth boundary. By taking the normal derivative of
Eq. \eqref{eq:IE} and letting the field point $\bm{x}$ going to the
boundary $\Gamma$, one obtains the hypersingular BIE (HBIE)
\begin{equation} \label{eq:HBIE}
    c(\bm{x})q(\bm{x})+\mathcal{H} u(\bm{x}) =\mathcal{M} q(\bm{x})+q^\text{inc}(\bm{x}), \quad \bm{x}\in \Gamma.
\end{equation}
Both CBIE and HBIE can be applied to solve the unknown
boundary values of interior problems. For exterior
problems, they have different set of fictitious frequencies
at which a unique solution can't be obtained.
However, Eqs. \eqref{eq:CBIE} and \eqref{eq:HBIE} will always
have only one solution in common. Given this fact, the Burton-Miller
formulation which is a linear combination of Eqs. \eqref{eq:CBIE}
and \eqref{eq:HBIE} (CHBIE) should yield unique solutions for
all frequencies
\begin{equation} \label{eq_bm}
\begin{split}
	c(\bm{x})u(\bm{x}) + \left(\mathcal{D} + \alpha\mathcal{H} \right)u(\bm{y}) - u^\text{inc}(\bm{x})
	=\left(\mathcal{S} + \alpha \mathcal{M} \right) q(\bm{y}) - \alpha \left[c(\bm{x}) q(\bm{x})
	- q^\text{inc}(\bm{x}) \right], \quad \bm{x}\in \Gamma,
\end{split}
\end{equation}
where, $\alpha$ is a coupling constant that can be chosen as $\text{i}/k$ \cite{BurtonMiller}.

\subsection{Nystr\"om boundary element discretization} \label{BEM:Ny-sin}

In this paper, the Nystr\"om method is used to discretize the Burton-Miller equation.
We choose to use the Nystr\"{o}m method because the resultant linear systems are
more like a summation, namely, the far-field part of the system matrix is
just the values of the integral kernels. As such, many fast methods, like
the FMM and the FDA, which are primarily proposed for
accelerating summations, can be applied to BEM acceleration with minor
modifications.

The implementation of the Nystr\"{o}m BEM follows the work in \cite{lcnystrom}. The boundary $\Gamma$ is partitioned into $N_\text{e}$ curved triangular quadratic elements. The 6-point Gauss quadrature rule on triangle is used
in evaluating regular element integrals, and thus the quadrature points of the Nystr\"{o}m method on each element are those of the  6-point Gauss rule. As a result, the total number of DOFs is $N=6\cdot N_\text{e}$.

For each quadrature point $\bm{x}_i, \, i=1,\cdots, \, N$, a local region, denoted by $D_i$, should be assigned. In this paper $D_i$ consists of all the elements whose distance to $\bm{x}_i$ is not larger than 2 times of their largest side length. To show the discretizing procedures, let $\mathcal{K}$ be one of the integral operators in
\eqref{eq:operator} and $K$ be the associated kernel function. For any field point $\bm{x}_i$ the boundary integral $\int_{\Gamma} K(\bm{x}_i, \bm{y}) u(\bm{y}) \mathrm{d}\bm{y}$ can be divided into two parts by the local region $D_i$; i.e.,
\begin{equation} \label{eq:ndbie}
    \int_{\Gamma} K(\bm{x}_i, \bm{y}) u(\bm{y}) \mathrm{d}\bm{y} \approx \sum_{\Delta \not\in D_i} \int_{\Delta} K(\bm{x}_i, \bm{y}) u(\bm{y}) \mathrm{d}\bm{y} + \sum_{\Delta \in D_i} \int_{\Delta} K(\bm{x}_i, \bm{y}) u(\bm{y}) \mathrm{d}\bm{y}\mathrm{d}\bm{y}
\end{equation}
where, $\Delta$ denotes the boundary elements. For elements outside the local region $D_i$ the integral is regular and thus is accurately evaluated by using the Gauss quadrature,
\begin{equation} \label{eq:far field}
    \int_{\Delta} K(\bm{x}_i,\bm{y})u(\bm{y}) \mathrm{d}\bm{y}
	\approx \sum_{j} \omega_j K(\bm{x}_i,\bm{y}_j) u(\bm{y}_j), \quad
	\Delta \not \in D_i,
\end{equation}
where, $\bm{y}_j$ and $\omega_j$ are the $j$-th quadrature point and weight
over element $\Delta$, respectively.

For elements inside the local region $D_i$, however, the kernels exhibit various types of
singularity. As a result, conventional
quadratures fail to give correct results. In order to maintain
high-order properties, the quadrature weights are adjusted
by a local correction procedure \cite{lcnystrom}. Consequently one has
\begin{equation} \label{eq:near field}
    \int_{\Delta} K(\bm{x}_i,\bm{y})u(\bm{y}) \mathrm{d} \bm{y}
	\approx \sum_j \bar\omega_j^{\Delta}(K; \bm{x}_i) u(\bm{y}_j), \quad \Delta \in D_i,
\end{equation}
where, $\bar\omega_j^{\Delta} (K; \bm{x}_i)$ represents the locally corrected quadrature
weights associated with element $\Delta$, kernel $K$ and field point $\bm{x}_i$. The local corrected procedure is performed by approximating the
unknown quantities using linear combination of polynomial basis functions which are defined on intrinsic coordinates of the element.
The locally corrected quadrature
weights are obtained by solving the linear system
\begin{equation} \label{eq:local correction}
    \sum_j \bar\omega_j^{\Delta} (K; \bm{x}_i) \phi^{(n)}(\bm{y}_j) = \int_\Delta
	K(\bm{x}_i,\bm{y})\phi^{(n)}(\bm{y})\mathrm{d} \bm{y}, \quad
	\Delta  \in D_i,
\end{equation}
where $\phi^{(n)}$ are polynomial basis functions. For the
Nystr\"om method based on quadratic elements as used in this paper, $\phi^{(n)}$ are given by
\begin{equation}
    \phi^{(n)}(\xi_1,\xi_2)=\xi_1^p\xi_2^q, \quad p+q \le 2,
\end{equation}
where, $p$ and $q$ are integers, $\xi_1$ and $\xi_2$ denote local intrinsic coordinates.

The integrals in Eq. \eqref{eq:local correction} become nearly singular when $\bm{x}_i$ is close to the element, and these integrals are computed by using a recursive subdivision quadrature in this paper. When the field point $\bm{x}_i$ locates on the element, singularity appears in the integrals. For the first three operators in equation \eqref{eq:operator}, the integrals have weak singularity of order $r^{-1}$,
while for the hyper-singular operator $\mathcal{H}$ the integral has singularity of order $r^{-3}$, as $r\to0$.
The accurate evaluation of those singular integrals is crucial to ensuring the accuracy of the BEM. Here, an efficient numerical method recently proposed in \cite{rjj_int} is used. The method is capable of treating weakly, strongly and hyper-singular integrals in a similar manner with high accuracy, and the code is open.

Consequently, the boundary integral \eqref{eq:ndbie} can be transformed into a summation as
\begin{equation} \label{eq:sum}
    \int_{\Gamma} K(\bm{x}_i, \bm{y}) u(\bm{y}) \mathrm{d}\bm{y} \approx \sum_{j=1}^{N} \bar K(\bm{x}_i, \bm{y}_j) u_j,
\end{equation}
where, $u_j=u(\bm{y}_j)$ and
\begin{equation} \label{kcorrection}
 \bar K(\bm{x}_i, \bm{y}_j) = \left\{ %
  \begin{array}{ll}
    K(\bm{x}_i, \bm{y}_j)\omega_j , & \bm{y}_j \hbox{ is on a element } \Delta \not\in D_i ,\\
    \bar\omega_j^{\Delta} (K; \bm{x}_i), & \bm{y}_j \hbox{ is on a element } \Delta \in D_i.
  \end{array}
\right.
\end{equation}

Performing the Nystr\"om discretization to the boundary integrals in the
Burton-Miller equation \eqref{eq_bm} leads to the BEM linear system of form
\begin{equation} \label{bemsys}
\bm{Hu=Gq+f},
\end{equation}
where, $\bm{H}$ and $\bm{G}$ are $N \times N$ matrices, $\bm{u}$ and $\bm{q}$ are $N$-vectors of the boundary values of $u(\bm{x})$ and $q(\bm{x})$, and $\bm{f}$ consists of the values of the incident wave terms in \eqref{eq_bm}. It is worth noting that when $\bm{y}_j$ is on a element outside the local region $D_i$ of $\bm{x}_i$ the kernel values
 of the matrices $\bm{H}$ and $\bm{G}$ are just the quadrature weight $\omega_j$ of $\bm{y}_j$ time the function values
\begin{subequations} \label{eq_hgmat}
\begin{align}
H_{ij} &= \frac{\partial G}{\partial \bm{n}_{\bm{y}}} (\bm{x}_i, \bm{y}_j) +
	\alpha \frac{\partial^2 G}{\partial \bm{n}_{\bm{x}} \partial \bm{n}_{\bm{y}}}
	(\bm{x}_i, \bm{y}_j) \quad \text{and}
	\label{eq_hmat} \\
G_{ij} &= G(\bm{x}_i, \bm{y}_j) +
	\alpha \frac{\partial G}{\partial \bm{n}_{\bm{x}}} (\bm{x}_i, \bm{y}_j).
	\label{eq_gmat}
\end{align}
\end{subequations}
Therefore, the matrix-vector product with matrices $\bm{H}$
and $\bm{G}$ is in fact summations with the respective kernels.

By considering the boundary conditions, equation \eqref{bemsys} is recast as the following linear system of equations to be solved
\begin{equation} \label{Axb}
\bm{Ax=b},
\end{equation}
where, the $N \times N$ matrix $\bm{A}$ consists of columns of matrices $\bm{H}$ and $\bm{G}$, depending on the specific boundary conditions, $\bm{b}$ and $\bm{x}$ are known and unknown $N$-dimensional vectors. For large-scale problems, linear system \eqref{Axb} is often solved by using the iterative solvers; the generalized minimal residual method (GMRES) will be used in this paper. The main computational work of the iterative solver is in the evaluation of matrix-vector product $\bm{Ax}$. Since matrix $\bm{A}$ is always fully populated in BEM, the computational cost for a naive evaluation should be $O(N^2)$, which suddenly becomes prohibitive with the increase of $N$. In the next section, the fast directional algorithm is employed to reduce the square scaled computational cost to almost linear.

\section{Fast directional algorithm for Burton-Miller formulation}
\label{s:fda}

In this section, the FDA is adapted to accelerate
the Nystr\"om BEM for Burton-Miller BIE, resulting in a fast BEM solver (denoted by FDBEM) for acoustic problems.
The central is to develop a FDA for the fast summation with the kernels in \eqref{eq_hgmat}.

\subsection{Multilevel FDA} \label{S-fda}

The multilevel FDA recently developed for the Helmholtz kernel is briefly reviewed here, detailed description
can be found in \cite{fda, fda2010}.

\subsubsection{Octree structure}

The implementation of the FDA relies on an adaptive octree structure
by which the source and target points are grouped \cite{kifmm}. 
The octree is constructed by the following steps. First, find a 
level-0 cube containing all the points. Then, subdivide each cube $C$ 
at level-$l$ of the octree into eight equally sized level-$(l+1)$ cubes
if it contains more than $N_\text{p}$ points. The subdivision process 
is performed recursively until each leaf cube contains no more than 
$N_\text{p}$ points. In this paper, $N_\text{p}$ is determined by
\begin{equation} \label{eq_npleaf}
	N_\text{p} = \max( 30, 50(\log_{10} \varepsilon - 3) ),
\end{equation}
where $\varepsilon$ is the target accuracy of the FDA.
The finest level of the octree is indexed by $L$.

In order to developing a FDA that is efficient in both low- and 
high-frequencies, the octree is often divided into two regimes, namely 
the low frequency regime and high frequency regime, according the side 
length of the cubes and the wavelength. The low frequency
regime consists of cubes whose side length $w$ is less
than the wavelength, and the high frequency regime consists of the rest 
cubes. Obviously, when the non-dimensional wave number $kD$ (where 
$D$ is the typical size of the boundary $\Gamma$) is sufficiently 
small, all the levels would lie in the low frequency regime, and 
the FDA would degenerate to the kernel independent FMM\cite{kifmm}.

\begin{figure}[h]
	\centering
	\includegraphics[width=0.5\textwidth]{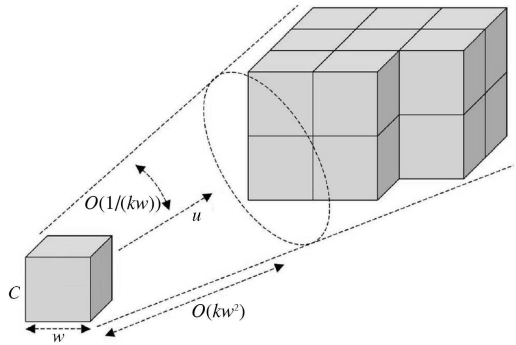}
	\caption{The interaction field
		of a cube $C$ in the high frequency regime is partitioned into multiple directional wedges.}
	\label{fig_para_sep_3d}
\end{figure}

With the aid of the octree, the basic notations of the FDA can be defined. 
Consider a cube $C$ with width $w$ in the octree.
\begin{itemize}
  \item \textit{Near field $N^C$:} When $C$ is in
the low frequency regime, its near field $N^C$ is defined as the union of
the cubes that are adjacent with $C$; when it is in the high frequency
regime, $N^C$ is defined as the union of the cubes $B$ that satisfy
$$
\text{dist}(B, C) \le c_{\text{d}}(k w^2),
$$
where, $\text{dist}(B, C)$ denotes the distance between $B$ and $C$; Let $w_B$ and $w_C$
	be the widths and $\bm{c}_B$ and $\bm{c}_C$ be the centroids of cubes $B$ and
	$C$, respectively, then
$$
\text{dist}(B, C)
	= \max_{i=1, 2, 3} (|\bm{c}_{B, i} - \bm{c}_{C, i}| - w_B/2 - w_C/2).$$
The constant $c_{\text{d}} = {1 \over \pi}$ is used in this paper.

  \item \textit{Far field $F^C$:} Consists of all the cubes not in $N^C$.

  \item \textit{Interaction field $I^C$:} $I^C = N^P \backslash N^C$.

  \item \textit{Directional wedges:} When $C$ is in
the high frequency regime, $I^C$
can be further divided into multiple directional wedges, each of which has the
spanning angle no greater than $O(1/(k w))$. The directional wedges would be
indexed by their center directions, as shown in Figure \ref{fig_para_sep_3d}. The low frequency regime can be considered as a special case of the high frequency regime which consists of only one
directional wedge indexed by $\bm{\zeta} = (0, 0, 0)$.

\end{itemize}
It is proved that the numerical rank of interaction matrices is low for the points in $C$ and each
directional wedge of $I^C$ \cite{fda}. Therefore, the far field summation can be accelerated by using the low rank approximation.

\subsubsection{Translations in FDA}

The FDA recently developed in \cite{fda} is a fast algorithm for computing summations with the Helmholtz kernel, i.e.,
\begin{equation}\label{eq:potential}
p_i = \sum_{j} G(\bm{x}_i, \bm{y}_j) q_j.
\end{equation}
where, $G$ is given by \eqref{eq:fundsol}. The basic procedure consists of five translations, namely, the source-to-multipole (S2M), multipole-to-multipole (M2M), multipole-to-local (M2L),
local-to-local (L2L) and local-to-target (L2T) translations, which is similar to the well-known FMM.
Here, the above five translations are briefly described.

\begin{figure}[h]
	\centering
	\subfigure[Low frequency regime.]{
		\includegraphics[width=0.35\textwidth]{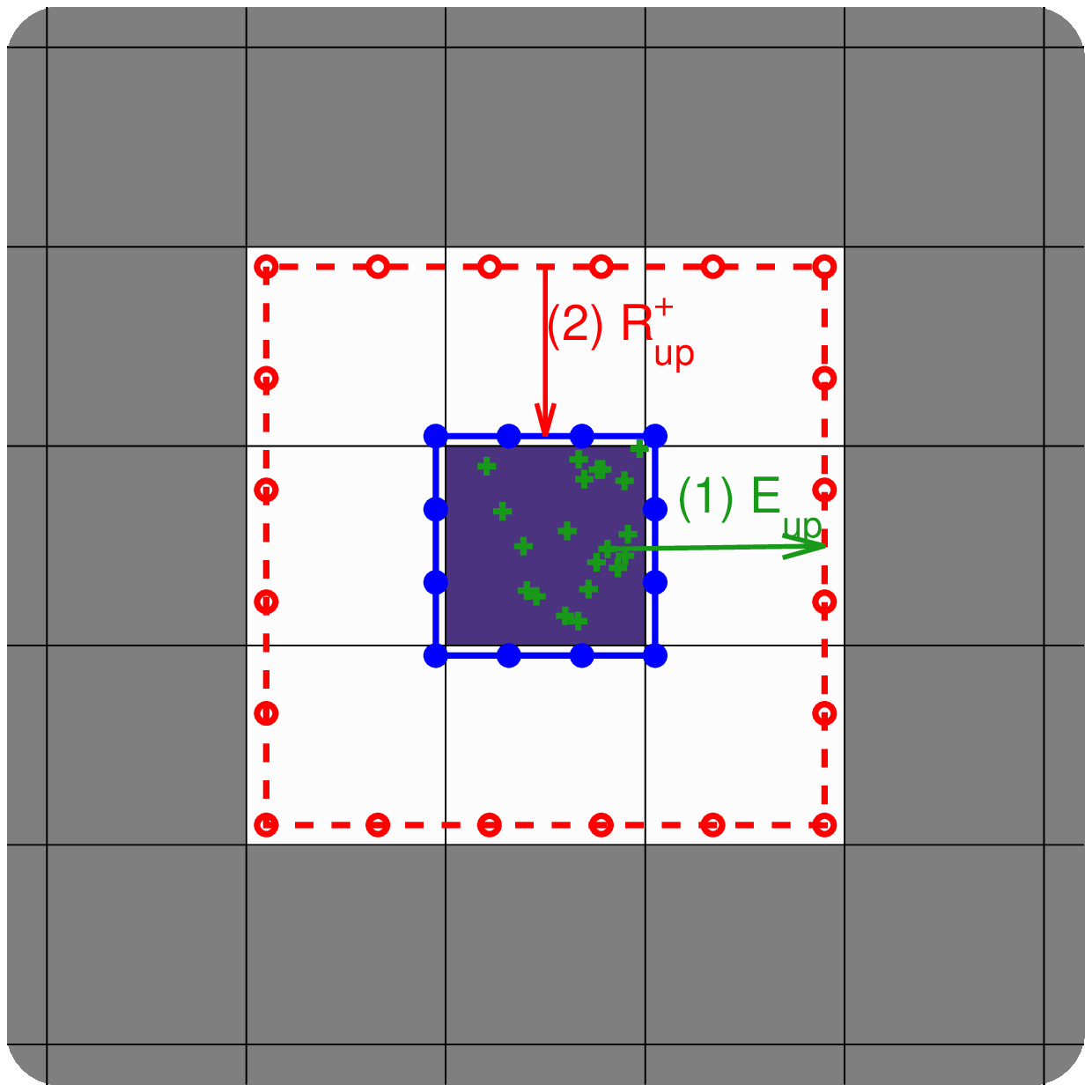}
		\label{fig_s2m_lfr}
	}
	\subfigure[High frequency regime.]{
		\includegraphics[width=0.5\textwidth]{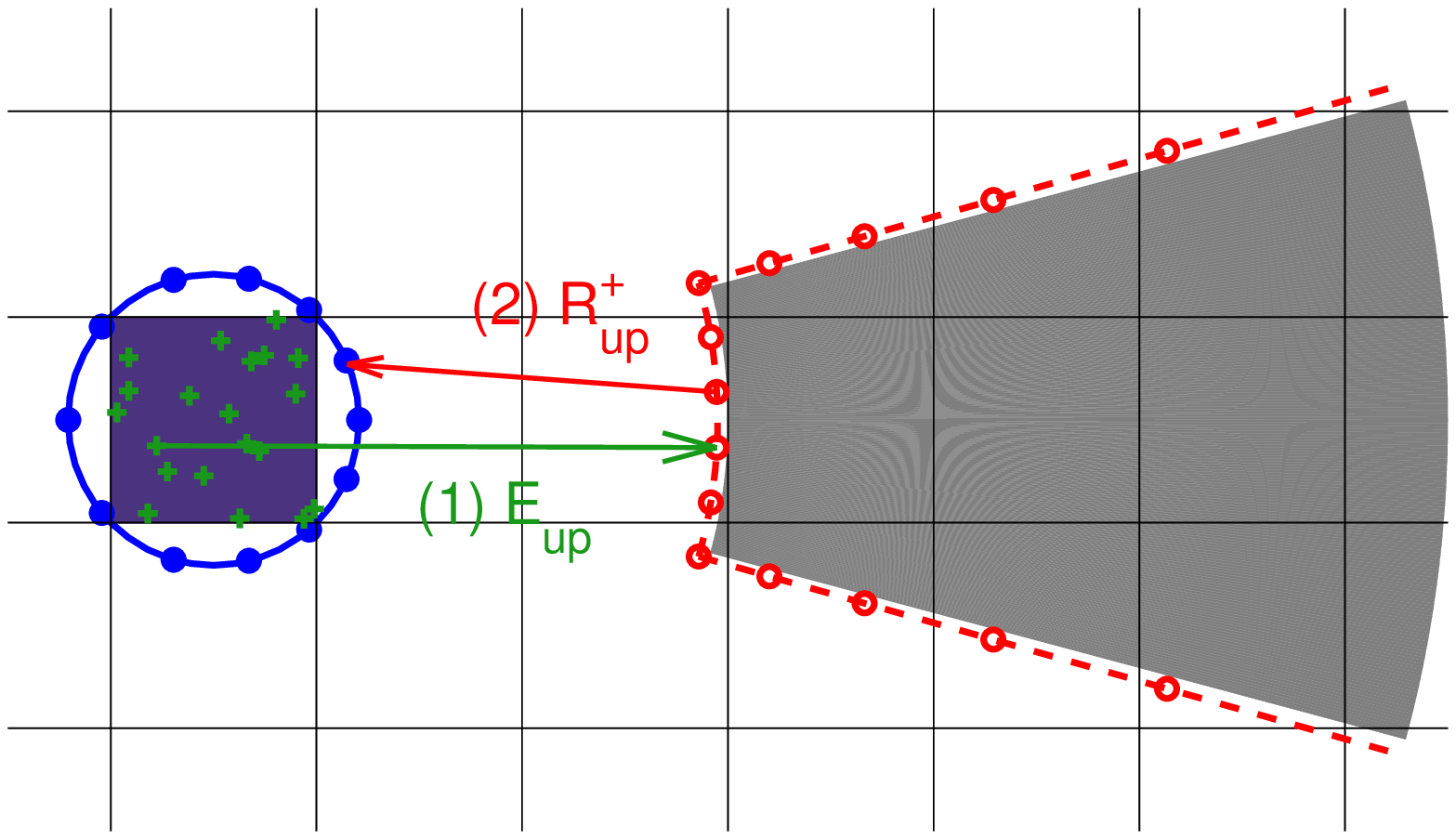}
		\label{fig_s2m_hfr}
	}
	\caption{The outgoing equivalent points (the blue points marked with
		``$\bullet$''), outgoing check points (the red points marked with
		``$\circ$'') and S2M translations in the low and high frequency regime.
		The computation of the outgoing equivalent densities for the original
		sources inside a cube (located at the points marked with ``+'') includes
		two steps shown by arrows: (1) the evaluation of the check potentials
		at the outgoing check points using the original sources, and (2) the
		inversion to construct the outgoing equivalent densities.}
	\label{fig_s2m}
\end{figure}

For a leaf cube $C$, the S2M operation, represented by matrix $\bm{S}$,
translates the original sources inside $C$ to the outgoing equivalent densities of $C$.
The definition of the outgoing equivalent points
and the associated outgoing check points are illustrated in Figure \ref{fig_s2m}, where
	all the definitions are demonstrated using the 2D cases, but they can be generalized to 3D readily. The S2M
translation can be further divided into two steps:
(1) the evaluation of the outgoing check potentials at the outgoing check points
produced by the original sources, the matrix of which is denoted by
$\bm{E}_\text{up}$, and (2) the inversion to construct the outgoing equivalent
densities at the outgoing equivalent
points that can reproduce the outgoing check potentials, the matrix of this linear translation is
denoted by $\bm{R}_\text{up}^+$.

When $C$ is in the low frequency regime, the outgoing equivalent points and
the outgoing check points are non-directional, and can be sampled straightforwardly
at Cartesian grid points, as described in \cite{kifmm}.
The matrix $\bm{R}_\text{up}^+$ is defined as the pseudo-inverse of
$[\bm{R}_\text{up}]_{ij} = G(\bm{a}_i, \bm{b}_j)$, where $\bm{a}_i$ is the $i$-th
outgoing check point, and $\bm{b}_j$ is the $j$-th outgoing equivalent point.
When $C$ is in the high frequency regime, since its interaction field is partitioned
into multiple directional wedges, a group of the outgoing equivalent and check points
have to be defined for each directional wedge.
In this paper, these points and $\bm{R}_\text{up}^+$ in direction $(1, 0, 0)$ are
computed by the algorithm in \cite{fda2010}. For the other directional wedges, the
points can be obtained by rotation and the matrices $\bm{R}_\text{up}^+$ remain the same.
In both low and high frequency regime, the numbers of equivalent
points and check points are controlled by the required accuracy
$\varepsilon$ of the FDA; see \cite{fda, fda2010}.

For each non-leaf cube at level-$l$, the M2M operation, associated with matrix $\bm{M}_l$, translates
the outgoing equivalent densities of its children to the outgoing equivalent densities of itself.
When $C$ is in the low frequency regime, only one group of outgoing equivalent densities of $C$
are to be translated, otherwise one has to loop over all the directional wedges
to compute the directional outgoing equivalent densities for each wedge.

\begin{figure}[h]
	\centering
	\subfigure[Low frequency regime.]{
		\includegraphics[width=0.35\textwidth]{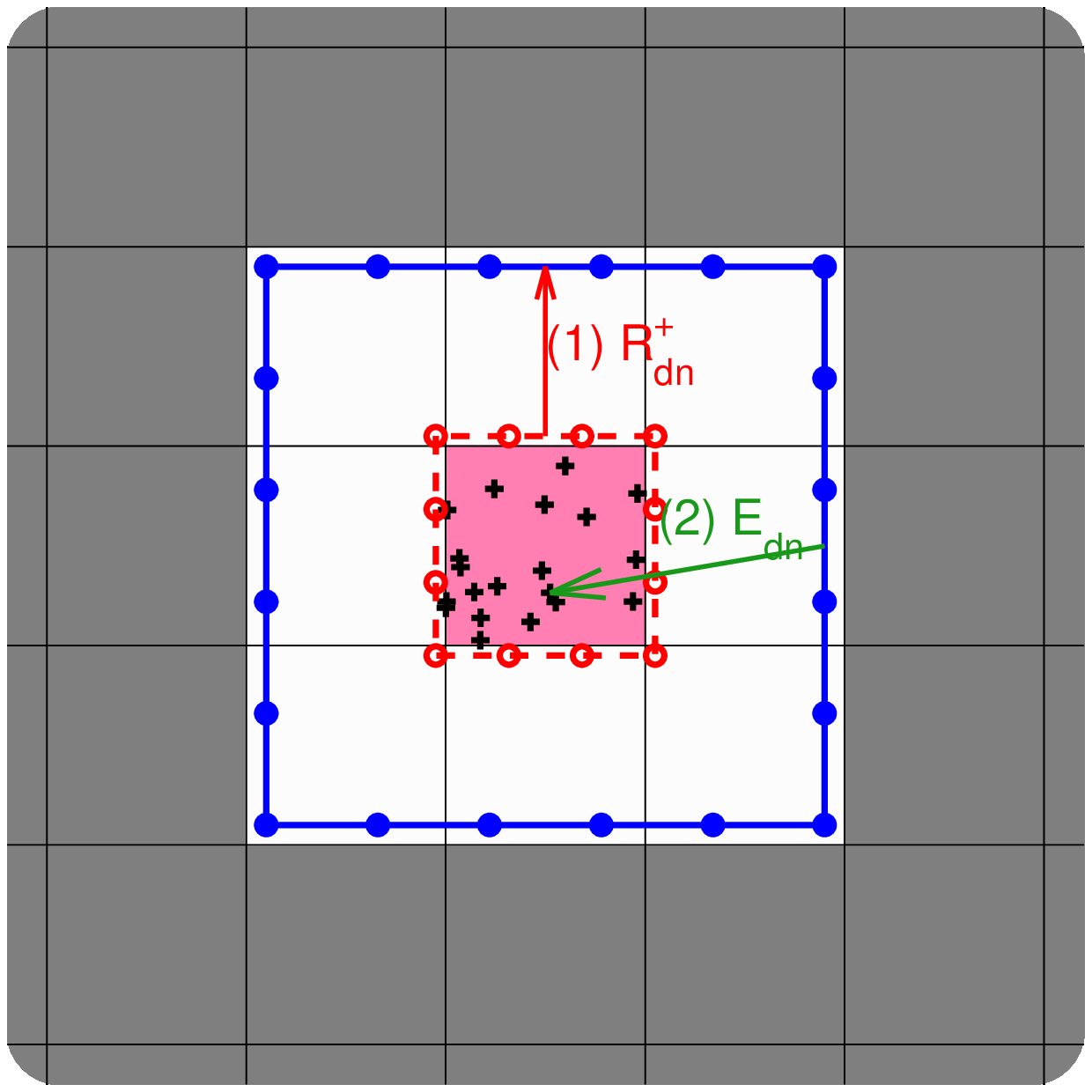}
		\label{fig_l2t_lfr}
	}
	\subfigure[High frequency regime.]{
		\includegraphics[width=0.5\textwidth]{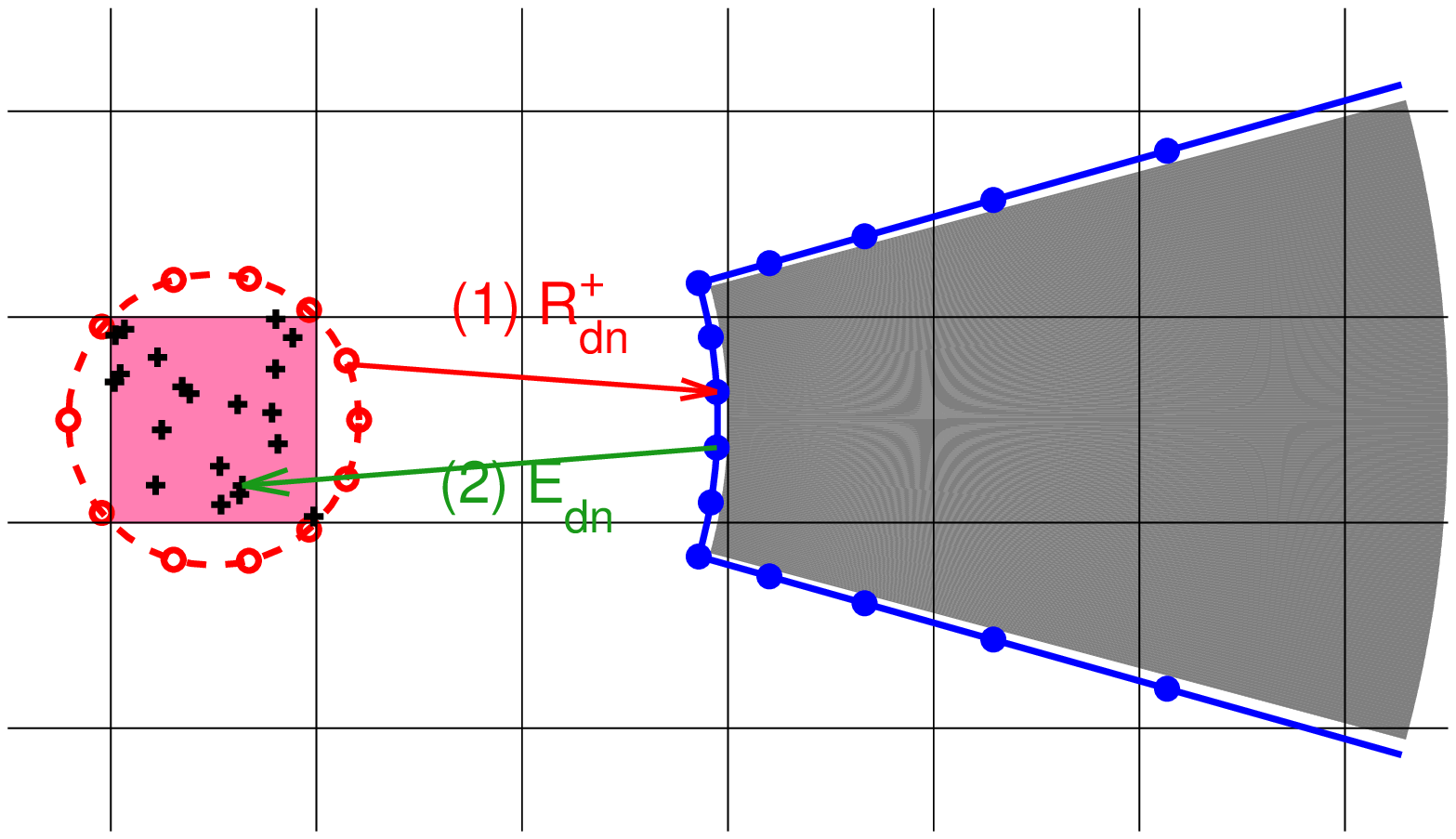}
		\label{fig_l2t_hfr}
	}
	\caption{The incoming equivalent points (the blue points marked with ``$\bullet$''), incoming
		check points (the red points marked with ``$\circ$'') and L2T translations in the low and high
		frequency regime. The potentials at the points inside the cube (marked with ``+'')
		is computed by two steps shown by arrows:
		(1) the inversion to construct the incoming equivalent densities that can reproduce the
		incoming check potentials at the check points, and (2) the evaluation of the potentials
		at the target points inside the cube using the equivalent densities.}
	\label{fig_l2t}
\end{figure}

The translations in the downward pass of FDA
are performed based on another two groups of points, namely
the incoming check points and incoming equivalent points whose definitions are similar to those of the outgoing points with the roles reversed. For example, in the downward pass the outgoing equivalent
points are used as the incoming check points, and the outgoing check points are used as
the incoming equivalent points, as illustrated in Figure \ref{fig_l2t}.

The M2L operation, represented by matrix $\bm{K}$, translate the outgoing equivalent points to incoming check potentials.
Assume $B$ is a cube in $C$'s interaction list and in direction $\bm{\zeta}$, then $C$ must lies in $B$'s
interaction list and in direction $-\bm{\zeta}$.
The M2L translation evaluate $B$'s incoming check potentials in direction $-\bm{\zeta}$
using $C$'s outgoing equivalent densities in direction $\bm{\zeta}$.

For each non-leaf cube $C$ at the $l$-th level, the L2L matrix $\bm{L}_l$ translates
the incoming check potentials of itself to the incoming check potentials of its children.
Similar to the translations in the upward pass, it
also consists of two steps: (1) the inversion to construct
$C$'s incoming equivalent densities at the incoming equivalent points that can reproduce
its incoming check potentials, the matrix of which is denoted by $\bm{R}_\text{dn}^+$, and
(2) the evaluation of the incoming check potentials of $C$'s
child cubes produced by $C$'s incoming equivalent densities, the matrix of which is denoted by $\bm{E}_\text{dn}$.

The L2T matrices $\bm{T}$ translate the incoming check potentials of each leaf cube $C$
to the potentials at the target points inside $C$. Similar to L2L, it also consists of
the inversion and evaluation steps, as illustrated in Figure \ref{fig_l2t}.

\subsection{FDA for Burton-Miller formulation}
\label{ss_modifications}

The FDA outlined above is a multilevel algorithm for the evaluation of potentials
defined by \eqref{eq:potential} where $G$ is the
Helmholtz kernel. Although the algorithm is kernel-independent, its application
to the evaluation of the potentials associated with the kernels \eqref{eq_hgmat}
of the Burton-Miller BIE is not straightforward, because the kernels involve
the normal vector $\bm{n_x}$ at the field point $\bm{x}$ which cannot be taken
into consideration in the construction of equivalent densities. Below, strategies
to cope with this issue and their numerical
implementations are described.

Consider the potential summation with kernel $G_{ij}$ in \eqref{eq_gmat},
\begin{equation} \label{eq_gsum}
\begin{split}
	p_i =& \sum_{j=1}^N \left[ G
	(\bm{x}_i, \bm{y}_j) + \alpha \frac{\partial G}{\partial \bm{n}_{\bm{x}}}
	(\bm{x}_i, \bm{y}_j) \right] q_j \\
	=& \sum_{j=1}^N \left( 1 + \alpha \frac{\partial}{\partial \bm{n_x}} \right)
	G (\bm{x}_i, \bm{y}_j) q_j \\
	=& \left( 1 + \alpha \frac{\partial}{\partial \bm{n_x}} \right) \sum_{j=1}^N
	G(\bm{x}_i, \bm{y}_j) q_j.
\end{split}
\end{equation}
It follows that the potential $p_i$ can be computed by two steps. The first
step is the evaluation of the single layer potentials, denoted by
$\hat p (\bm{x})$, generated by the sources $q_j$, i.e.,
\begin{equation}	\label{eq_sump0g}
	\hat p(\bm{x}) = \sum_{j=1}^N G(\bm{x}, \bm{y}_j) q_j.
\end{equation}
This can be accelerated by using the FDA in Section \ref{S-fda}.
The second step is the computation of $p_i$ by
\begin{equation} \label{eq_pip0}
	p_i = \left( 1 + \alpha \frac{\partial}{\partial \bm{n_x}} \right)
	\hat p(\bm{x}_i).
\end{equation}
The numerical implementation of this step requires a change of the matrix
$\bm{E}_\text{dn}$ that is used in the computation of the L2T translation matrix
$\bm{L} = \bm{E}_\text{dn} \bm{R}^+$. More specifically, for a cube $C$, suppose
that a group of incoming equivalent densities $\bar{q}_j$ at the incoming
equivalent points $\bm{b}_j$ are obtained in the first step of the L2T
translation, which can produce $\hat p(\bm{x})$ inside $C$, i.e.,
$$
\hat{p}(\bm{x}) = \sum_j G(\bm{x}, \bm{b}_j) \bar{q}_j.
$$
Then the potential $p_i$ in \eqref{eq_pip0} can be evaluated by
\begin{equation}
\begin{split}
	p_i =& \left( 1 + \alpha \frac{\partial}{\partial \bm{n_x}} \right)
	\sum_j G(\bm{x}_i, \bm{b}_j) \bar{q}_j \\
	=& \sum_j \left( 1 + \alpha \frac{\partial}{\partial \bm{n_x}} \right)
	G(\bm{x}_i, \bm{b}_j) \bar{q}_j \\
	=& \sum_j \left[ G(\bm{x}_i, \bm{b}_j) + \alpha \frac{\partial G}
	{\partial \bm{n_x}}(\bm{x}_i, \bm{b}_j) \right] \bar{q}_j \\
\end{split}
\end{equation}
This suggests that $p_i$ at $\bm{x}_i$ can be computed by using
\begin{equation} \label{eq_eval2t}
	[\bm{E}_\text{dn}]_{ij} = G(\bm{x}_i, \bm{b}_j) + \alpha \frac{\partial G}
	{\partial \bm{n_x}}(\bm{x}_i, \bm{b}_j)
\end{equation}
as the translation matrix for the evaluation step in L2T.

The FDA for the potential summation with kernel $H_{ij}$ in \eqref{eq_hmat} can be attained analogously. Since
\begin{equation} \label{eq_hsum}
\begin{split}
	p_i =& \sum_{j=1}^N \left[ \frac{\partial G}{\partial \bm{n}_{\bm{y}}}
	(\bm{x}_i, \bm{y}_j) + \alpha \frac{\partial^2 G}{\partial \bm{n}_{\bm{x}}
	\partial \bm{n}_{\bm{y}}} (\bm{x}_i, \bm{y}_j) \right] q_j \\
	=& \sum_{j=1}^N \left( 1 + \alpha \frac{\partial}{\partial \bm{n_x}} \right)
	\frac{\partial G}{\partial \bm{n_y}} (\bm{x}_i, \bm{y}_j) q_j \\
	=& \left( 1 + \alpha \frac{\partial}{\partial \bm{n_x}} \right) \sum_{j=1}^N
	\frac{\partial G}{\partial \bm{n_y}} (\bm{x}_i, \bm{y}_j) q_j.
\end{split}
\end{equation}
The potential $p_i$ can be computed with
two steps. The first step is
the evaluation of the double layer potentials, denoted by $\tilde p (\bm{x})$,
generated by the sources $q_j$,
\begin{equation}	\label{eq_sump0h}
	\tilde p(\bm{x}) = \sum_{j=1}^N \frac{\partial G}{\partial \bm{n_y}}
	(\bm{x}, \bm{y}_j) q_j.
\end{equation}
Summation \eqref{eq_sump0h} can be accelerated by using the FDA in Section \ref{S-fda} with minor modification. The
rationale is that the source densities $q_j$ are essentially dipoles in the
$\bm{n_y}$ direction. Since the field produced by dipoles can be approximated by
a group of monopoles in its vicinity, the potentials produced by the
sources inside a cube $C$ can be approximated by a group of equivalent
densities, when observing at the far field of $C$.

Therefore, in the FDA for single layer
kernel $G$ in Section \ref{S-fda}, the S2M matrix $\bm{S}$ is computed as
$\bm{S} = \bm{R}_\text{up}^+ \bm{E}_\text{up}$. The entries of matrix
$\bm{E}_\text{up}$ are given by the values of the single layer kernel,
$$
[\bm{E}_\text{up}]_{ij} = G(\bm{a}_i, \bm{y}_j),
$$
because the source densities are monopoles. When FDA is adopted to the double
layer potential in \eqref{eq_sump0h}, the entries of $\bm{E}_\text{up}$ have
to be the values of the double layer kernel,
\begin{equation}	\label{eq_evas2mh}
	[E_\text{up}]_{ij} = \frac{\partial G}{\partial \bm{n_y}} (\bm{a}_i, \bm{y}_j).
\end{equation}

The second step is the computation of $p_i$ by
\begin{equation}
	p_i = \left( 1 + \alpha \frac{\partial}{\partial \bm{n_x}} \right)
	\tilde p(\bm{x}_i),
\end{equation}
which can be realized in the same way as \eqref{eq_pip0}.

From the above discussion, it is clear that there are two common points in the FDA
for the Helmholtz kernel \eqref{eq:potential} and those for the kernels of the
Burton-Miller BIE. First, equivalent sources of monopoles are always used even
though the Burton-Miller kernels \eqref{eq_hgmat} are the combinations of the
Helmholtz kernel and its derivatives. Second, since the M2M, M2L and L2L
translations only involve the equivalent densities, no modification is needed
at all for these translations.

Besides the common points, three modifications are described for the application
of FDA to the Burton-Miller formulation:
\begin{enumerate}
	\item[(1)] For the summation with $H_{ij}$, the outgoing check potentials for
	S2M should be computed by using the kernel
	$\frac{\partial G}{\partial \bm{n_y}} (\bm{x}, \bm{y})$.
	\item[(2)] The matrix $\bm{E}_\text{dn}$ in the computation of the L2T matrix
	should be evaluated by \eqref{eq_eval2t}, for both summations with kernels
	$H_{ij}$ and $G_{ij}$.
	\item[(3)] The contribution to the potentials by the sources in the near
	field should be computed by using the kernels \eqref{eq_hmat} and
	\eqref{eq_gmat}, as well as their local corrections.
\end{enumerate}

In the implementation of the FDA, the matrices of the M2M, M2L and L2L translations
are computed only once and then stored in memory for later use. In addition, from
the definition of the equivalent points and check points, it is known that
the inverse operator in the downward pass is the transpose of that in
the upward pass, and the the L2L translation matrix is just the transpose of
the M2M translation matrix for the same directional wedge; that is,
$\bm{R}^+_{\text{dn}} = (\bm{R}^+_{\text{up}})^\text{T}$, $\bm{E}_\text{dn}
= \bm{E}_\text{up}^\text{T}$, and
\begin{equation}
\bm{L} = \bm{E}_\text{dn} \bm{R}^+_\text{dn}
= \bm{E}_\text{up}^\text{T} (\bm{R}^+_\text{up})^\text{T}
= (\bm{R}^+_\text{up} \bm{E}_\text{up})^\text{T}
= \bm{M}^\text{T}.
\end{equation}
Therefore only one set of matrix is needed to be computed and saved for both
the L2L and M2M translations in each directional wedge.

\section{Further improvement to FDA}
\label{s:improve}

\begin{figure}[h]
	\centering
	\includegraphics[width=0.6\textwidth]{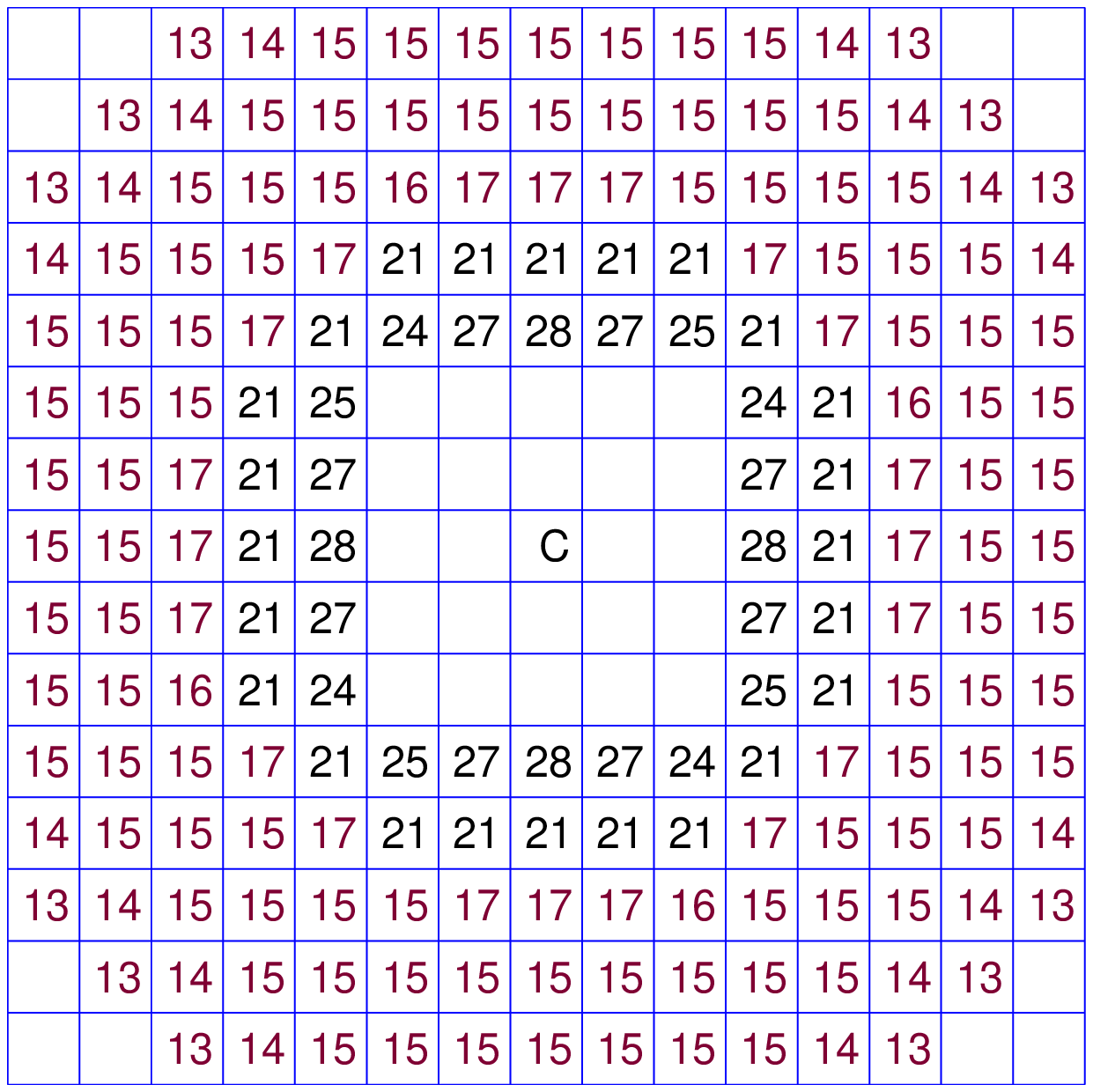}
	\caption{The ranks of M2L matrices in the finest level of the high frequency
		regime, i.e., $w=\lambda$, with $\varepsilon =$ 1e-4. The dimension
		of the compressed M2L matrices is $s = 38$. The purple number
		denotes the M2L matrices whose rank is less than $s/2$, that is, these
		translations can be accelerated by using the low rank approximations.}
	\label{fig_rankm2l}
\end{figure}

In this section, further improvement to the FDA is proposed based on the observation
that generally the translation matrices are rank deficient. Therefore,
the computational cost of the FDA can be reduced by compressing these translation matrices
into more compact forms. Furthermore, our numerical experiments show that, after
the matrix reduction, most of the M2L matrices are still of ranks much lower than
their dimensions, as shown in Figure \ref{fig_rankm2l}. Therefore these M2L
translations can be further accelerated by using the low rank approximations for
the translation matrices. In the present work, due to the small sizes of the
matrices, the low-rank decomposition is computed by using singular value
decomposition (SVD) to get the optimal ranks.

\subsection{Matrix reduction}
\label{sbsbsec_matred}

Consider two leaf cubes $B$ and $C$ at the $l$-th level such that $B$ is in the interaction field of $C$'s parent cube.
Then, the potentials
$\bm{p}_l$ in $C$ produced by the source densities $\bm{q}_l$ in $B$ are
evaluated as
\begin{equation} \label{eq_trans}
	\bm{p}_l = \bm{T}_{l} \bm{L}_{l-1} \bm{K}_{l-1} \bm{M}_{l-1} \bm{S}_{l} \bm{q}_l,
\end{equation}
where, $\bm{S}_l$ and $\bm{T}_l$ are the S2M and L2T translation matrices at the $l$-th level; $\bm{M}_{l-1}$, $\bm{K}_{l-1}$ and $\bm{L}_{l-1}$ are the M2M, M2L and L2L translation matrices at the $(l-1)$-th level.

Below, a scheme for reducing the sizes of aforementioned matrices is proposed.
Consider the computation of the S2M and L2T matrices,
\begin{subequations} \label{eq_matst}
\begin{equation}
	\bm{S}_l = \bm{R}_{\text{up}, l}^{+} \bm{E}_{\text{up}, l},
\end{equation}
\begin{equation}
	\bm{T}_l = \bm{E}_{\text{dn}, l} \bm{R}_{\text{dn}, l}^{+},
\end{equation}
\end{subequations}
where, $\bm{R}^+$ is the pseudo-inverse of $\bm{R}$, which is computed by the SVD,
\begin{equation} \label{eq_svdr}
	\bm{R} = \bm{U} \bm{\Sigma} \bm{V}^\text{H},
\end{equation}
where, $\bm{U}$ and $\bm{V}$ are unity matrices, $\bm{V}^\text{H}$
is the Hermitian matrix of $\bm{V}$, and $\bm{\Sigma} = \text{diag}(\sigma_1,
\sigma_2, \cdots, \sigma_n),\,\text{with} \, \sigma_1 > \sigma_2 > \cdots >
\sigma_n$, is the diagonal matrix of singular values. For a given error tolerance
$\varepsilon$, identify the singular values $\sigma_i < \varepsilon \sigma_1$
and truncate the associated columns of matrices $\bm{U}$ and $\bm{V}$ to obtain
two matrices $\tilde{\bm{U}}$ and $\tilde{\bm{V}}$ which consist of the columns
of $\bm{U}$ and $\bm{V}$ corresponding to the singular values $\sigma_i \ge
\varepsilon \sigma_1$. Then $\bm{R}^+$ is given by
\begin{equation} \label{eq_invr}
	\bm{R}^{+} \approx \tilde{\bm{V}} \tilde{\bm{\Sigma}}^{-1}
	\tilde{\bm{U}}^\text{H},
\end{equation}
where,  $\tilde{\bm{\Sigma}}$ is the diagonal matrix of the remaining singular
values. In this paper, the error tolerance $\varepsilon$ is chosen to be the
controlling accuracy in the sampling process of the directional equivalent
points and check points.

The matrices $\bm{M}_{l-1}$ and $\bm{L}_{l-1}$ have the similar decompositions
with the matrices $\bm{S}_l$ and $\bm{T}_l$ in \eqref{eq_matst}. By plugging
those decompositions into \eqref{eq_trans} and using \eqref{eq_invr}, one gets
\begin{equation} \label{eq_cmptrans}
\begin{split}
	\bm{p}_l = & \bm{T}_{l} \bm{L}_{l-1} \bm{K}_{l-1} \bm{M}_{l-1} \bm{S}_{l} \bm{q}_l	\\
		= & (\bm{E}_{\text{dn}, l} \bm{R}_{\text{dn}, l}^+) (\bm{E}_{\text{dn}, l-1} \bm{R}_{\text{dn}, l-1}^+)
			\bm{K}_{l-1} (\bm{R}_{\text{up}, l-1}^+ \bm{E}_{\text{up}, l-1})
			(\bm{R}_{\text{up}, l}^+ \bm{E}_{\text{up}, l}) \bm{q}_l	\\
		\approx & \bm{E}_{\text{dn}, l} (\tilde{\bm{V}}_{\text{dn}, l} \tilde{\bm{\Sigma}}_{\text{dn}, l}^{-1} \tilde{\bm{U}}_{\text{dn}, l}^\text{H})
			\bm{E}_{\text{dn}, l-1} (\tilde{\bm{V}}_{\text{dn}, l-1} \tilde{\bm{\Sigma}}_{\text{dn}, l-1}^{-1} \tilde{\bm{U}}_{\text{dn}, l-1}^\text{H}) \\
		 &	\bm{K}_{l-1}
			(\tilde{\bm{V}}_{\text{up}, l-1} \tilde{\bm{\Sigma}}_{\text{up}, l-1}^{-1} \tilde{\bm{U}}_{\text{up}, l-1}^\text{H}) \bm{E}_{\text{up}, l-1}
			(\tilde{\bm{V}}_{\text{up}, l} \tilde{\bm{\Sigma}}_{\text{up}, l}^{-1} \tilde{\bm{U}}_{\text{up}, l}^\text{H}) \bm{E}_{\text{up}, l} \bm{q}_l	\\
		= & \tilde{\bm{T}}_l \tilde{\bm{L}}_{l-1} \tilde{\bm{K}} \tilde{\bm{M}}_{l-1} \tilde{\bm{S}}_l \bm{q}_l,
\end{split}
\end{equation}
where, all the new translation matrices have smaller sizes,
\begin{equation} \label{eq_cmps}
	\tilde{\bm{S}}_l = \tilde{\bm{\Sigma}}_{\text{up}, l}^{-1} \tilde{\bm{U}}_{\text{up}, l}^\text{H} \bm{E}_{\text{up}, l},
\end{equation}
\begin{equation} \label{eq_cmpm}
	\tilde{\bm{M}}_l = \tilde{\bm{\Sigma}}_{\text{up}, l}^{-1} \tilde{\bm{U}}_{\text{up}, l}^\text{H} \bm{E}_{\text{up}, l} , \tilde{\bm{V}}_{\text{up}, l+1}
\end{equation}
\begin{equation} \label{eq_cmpk}
	\tilde{\bm{K}}_l = \tilde{\bm{U}}_{\text{dn}, l}^\text{H} \bm{K}_{l} \tilde{\bm{V}}_{\text{up}, l},
\end{equation}
\begin{equation} \label{eq_cmpl}
	\tilde{\bm{L}}_l = \tilde{\bm{U}}_{\text{dn}, l+1}^\text{H} \bm{E}_{\text{dn}, l} \tilde{\bm{V}}_{\text{dn}, l} \tilde{\bm{\Sigma}}_{\text{dn}, l}^{-1}
	= \tilde{\bm{M}}_l^\text{T},
\end{equation}
\begin{equation} \label{eq_cmpt}
	\tilde{\bm{T}}_l = \bm{E}_{\text{dn}, l} \tilde{\bm{V}}_{\text{dn}, l} \tilde{\bm{\Sigma}}_{\text{dn}, l}^{-1}.
\end{equation}
Thus the computational cost of all the translations can be reduced by using these new matrices \eqref{eq_cmptrans}. It is shown that the reduction of computational time can be up to $30\%$ for each FDA summation. Below, a method for further and much more considerable reduction of the computational time of the FDA is proposed.

\subsection{Low rank approximations for the compressed M2L matrices}

Here the M2L translation is considered in detail because it consumes
the most computational time of the FDA. Let $\tilde{\bm{K}}$ be a
compressed M2L matrix in \eqref{eq_cmpk} with dimension $s$. It is
found that, to the accuracy of $\varepsilon$, there
are always many M2L matrices $\tilde{\bm{K}}$ whose numerical ranks
are much less than $s/2$. This indicates that the M2L translation would be further accelerated
by performing matrix-vector production using the low rank
decomposition of $\tilde{\bm{K}}$ instead of $\tilde{\bm{K}}$ itself.

The low rank decomposition for each compressed M2L matrix $\tilde{\bm{K}}$
is computed by the truncated SVD with the singular values $\sigma_i  < \varepsilon \sigma_1 $ being dropped; that is,
\begin{equation} \label{eq_lrdk}
	\tilde{\bm{K}} \approx \hat{\bm{U}} \hat{\bm{S}}
	\hat{\bm{Q}}^\text{H}
	= \hat{\bm{U}} \hat{\bm{V}},
\end{equation}
where, $\hat{\bm{U}}$ and $\hat{\bm{Q}}$ are matrices with orthogonal columns, $\hat{\bm{S}}$ are the truncated singular value matrix and $\hat{\bm{V}} = \hat{\bm{S}}\hat{\bm{Q}}^\text{H}$. Moreover, since there are at most $O(\sqrt{N})$ M2L matrices
in FDA \cite{fda}, the computational time of the decomposition of the M2L matrices is of order $O(s^3 \sqrt{N})$
which is negligible in the total CPU time of the FDA.

For each M2L matrix $\tilde{\bm{K}}$ with rank $r<s/2$, the low-rank decomposition \eqref{eq_lrdk} is used in the translation. Otherwise, the matrix $\tilde{\bm{K}}$ itself is used.

The objective of the compression and decomposition to the translation matrices in this section is similar with the SArcmp M2L optimizer in \cite{opm2l}. However, two distinctions of our method should be noticed.
First, the compressed matrices in section \ref{sbsbsec_matred}
are computed during the definition of the equivalent and check points, thus one do not need to compute the low rank decomposition of the collected M2L matrices as in SArcmp.
Second, all the translation matrices are compressed in our work, while only the M2L matrices are
compressed in \cite{opm2l}.

\section{Numerical studies} \label{s:ne}

In this section, representative numerical examples are provided to demonstrate
the performance of the improved FDA and its application for the acceleration of the BEM for Burton-Miller formulation, denoted by FDBEM. Our codes are implemented serially in C++.
The computing platform is a workstation with a Xeon 5450 (2.66 GHz) CPU and 32 GB RAM.

\subsection{Performance of the improvements in section \ref{s:improve}}
\label{sbsec_numimp}

First, the performance of our acceleration technique in
Section \ref{s:improve} is tested by evaluating the potential summation
\begin{equation} \label{eq_sum}
	p_i = \sum_{j = 1}^{N} \left[ G(\bm{x}_i, \bm{y}_j) + \alpha
	\frac{\partial G}{\partial \bm{n_y}} (\bm{x}_i, \bm{y}_j) \right] q_j,
\end{equation}
with $G$ being the Helmholtz kernel in \eqref{eq:fundsol}.
We set $\alpha = \text{i}/k$ as this will be used in the Burton-Miller formulation.
The points $\{\bm{x}_i\}$ and $\{\bm{y}_j\}$
are sampled on the surface of a unit sphere with about 20 points per wavelength.
The densities $\{q_j\}$ on $\{\bm{y}_j\}$ are randomly defined
 with mean 0.
Equation \eqref{eq_sum} is evaluated by the original FDA and the improved FDA in section \ref{s:improve}.
The error $\varepsilon_\text{a}$ of the potential $p_i$,
\begin{equation}
	\varepsilon_\text{a} = \sqrt{ \frac{\sum_{i=1}^{N_\text{t}} \left|
	p_i^\text{(a)} - p_i^\text{(d)} \right|^2}{\sum_{i=1}^{N_\text{t}}
		\left| p_i^\text{(d)} \right|^2} },
\end{equation}
is computed and compared, where $p_i^\text{(a)}$ is the potential on a $N_\text{t}$ randomly
selected points $\bm{x}_i$ computed by the fast algorithm, while
$p_i^\text{(d)}$ is the potentials evaluated by straightforwardly
summation. In this paper, $N_\text{t}$ is chosen as 200.
For particle summation problems, the S2M and L2T matrices are used
only once,
thus they are not pre-computed but only computed when they are to be used.

\begin{table} [h]
	\centering
	\caption{Results of the example in section \ref{sbsec_numimp}.}
	\vspace{1em}
	\begin{tabular}{rr|rr|rr|rr|rr}
		\hline
		\multirow{2}{*}{$N$}	& \multirow{2}{*}{$k$}	& \multicolumn{2}{|c}{$T_\text{M2L}$ (s)}	& \multicolumn{2}{|c}{$T_\text{up}$ (s)}	& \multicolumn{2}{|c}{ $\varepsilon_\text{a}$ }	& \multicolumn{2}{|c}{$M$ (MB)}	\\
		\cline{3-10}
				& 			& original	& improved	& original	& improved		& original	& improved		& original	& improved	\\
		\hline
		\multicolumn{10}{c}{$\varepsilon =$ 1e-4}	\\
		73728	& $8\pi$	&   10	&   4		&   5		&    3		& 4.86e-5	& 5.88e-5		&   310	&   158	\\
		294912	& $16\pi$	&   44	&  17		&  17		&   10		& 7.30e-5	& 1.31e-4		&   735	&   441	\\
		1179648	& $32\pi$	&  218	&  94		& 103		&   58		& 1.15e-4	& 2.33e-4		&  1993	&  1338	\\
		4718592	& $64\pi$	&  966	& 383		& 458		&  214		& 2.00e-4	& 3.70e-4		&  6818	&  5053	\\
		\multicolumn{10}{c}{$\varepsilon =$ 1e-6} 	\\
		73728	& $8\pi$	&   33	&   11		&   12		&    7		& 2.97e-7	& 9.33e-7		&  1010	&   471	\\
		294912	& $16\pi$	&  154	&   55		&   59		&   36		& 7.75e-7	& 2.72e-6		&  2186	&  1189	\\
		1179648	& $32\pi$	&  631	&  235		&  218		&  145		& 1.30e-6	& 3.83e-6		&  4744	&  2941	\\
		4718592	& $64\pi$	& 2937	& 1145		& 1153		&  742		& 3.26e-6	& 1.06e-5		& 13938	&  9764	\\
		\multicolumn{10}{c}{$\varepsilon =$ 1e-8} 	\\
		73728	& $8\pi$	&   93	&   29		&   24		&   16		& 2.39e-8	& 2.79e-8		&  2644	&  1156	\\
		294912	& $16\pi$	&  421	&  137		&  119		&   84		& 4.93e-8	& 1.11e-7		&  5154	&  2684	\\
		1179648	& $32\pi$	& 1806	&  602		&  563		&  384		& 7.12e-8	& 7.27e-8		& 10147	&  5944	\\
		4718592	& $64\pi$	& 7721	& 2527		& 2478		& 1586		& 1.71e-7	& 2.55e-7		& 26860	& 18121	\\
		\hline
	\end{tabular}
	\label{tab_improve}
\end{table}

The CPU times and errors are given in Table \ref{tab_improve}. $T_\text{M2L}$ and $T_\text{up}$ are the CPU time of the M2L translation and upward pass, respectively; $M$ is the total memory usage. The summations
are evaluated with different controlling accuracies $\varepsilon$ and wave number $k$.
It can be seen that for arbitrary frequencies the M2L translation of the improved
FDA is always 2 to 3 times faster than that of the original FDA.
The upward passes are accelerated by a factor of about 40\%. Since the translation
matrices in the downward pass have the same size with those in the upward pass, the
same improvement can be expected for the downward pass.
The improved FDA can reduce the memory consumption by about 30\%.
The resulting errors are increased slightly, since additional errors are
introduced in the acceleration technique.

\subsection{Efficiency and accuracy of the FDBEM}
\label{ss_neFDBEM}

Here the efficiency and accuracy of the FDBEM developed in this paper are
demonstrated by solving a benchmark problem: the sound radiation of a unit sphere pulsating with radial velocity
$q = \frac{\partial u}{\partial \bm{n}} = 1$. This is a Neumann boundary value problem. The velocity potential $u$ on the spherical surface has analytical solution $u = \frac{1}{1 - \text{i}k}$.

\begin{figure}[h]
	\centering
	\includegraphics[width=0.47\textwidth]{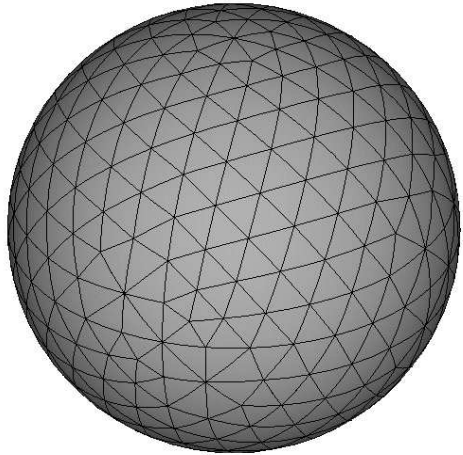}
	\caption{The spherical surface is discretized with 720 triangular
	quadratic elements for $k = 2\pi$.}
	\label{fig_usphmesh}
\end{figure}

Radiation problems with wave numbers $k = 2\pi, 4\pi, 8\pi, 16\pi$
and $32\pi$ are solved. For each problem, the spherical surface
is discretized by curved triangular quadratic elements with element
size $h = \lambda/5$, where $\lambda = 2\pi/k$ is the
wave length. In the problem with $k=2\pi$, the wave length is
one half of the diameter of the sphere. The mesh consists
of $N_\text{e} = 720$ elements (see Fig. \ref{fig_usphmesh}) and the DOF
is $N = 6 N_\text{e} = 4320$. In the problem with the highest wave
number $k=32\pi$, the wave length is 1/32 of the diameter. The
corresponding mesh consists of $N_\text{e} = 188796$ elements,
$N = 1132776$ DOFs.

One should notice that this series of wave numbers are exactly the characteristic frequencies of the unit sphere, at which the BEM based on the Helmholtz BIE fails to get the correct solution. However, as it will be shown that our FDBEM can still achieve highly accurate results.
In the FDBEM, the S2M and L2T matrices are
precomputed and saved in memory.
The resulting linear system is solved by the GMRES solver with
convergence tolerance being equal to $\varepsilon$.

\begin{figure}[h]
	\centering
	\subfigure[Relative error versus the DOF.]{
		\includegraphics[width=0.47\textwidth]{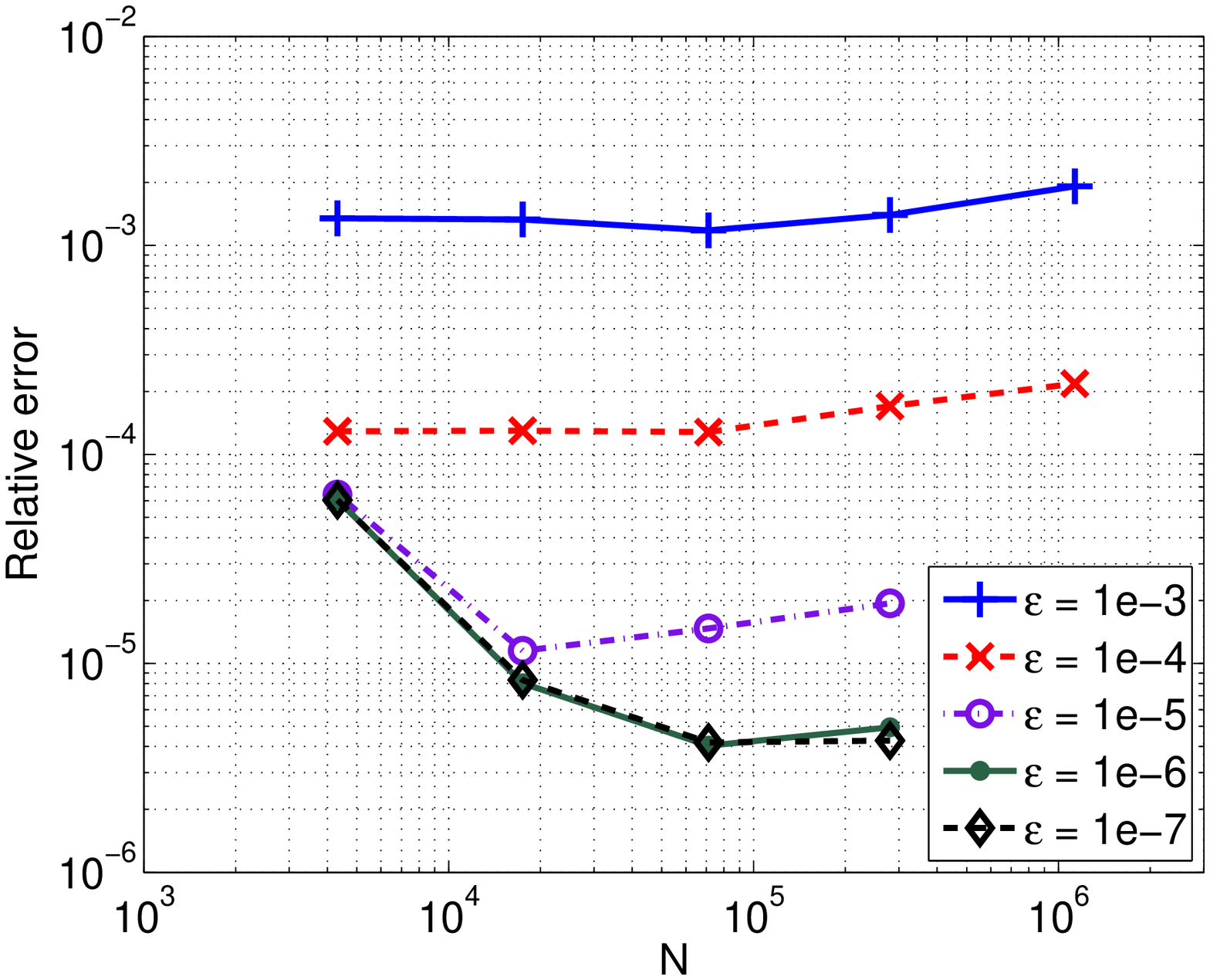}
		\label{fig_errbm}
	}
	\subfigure[Memory consumption versus the DOF.]{
		\includegraphics[width=0.47\textwidth]{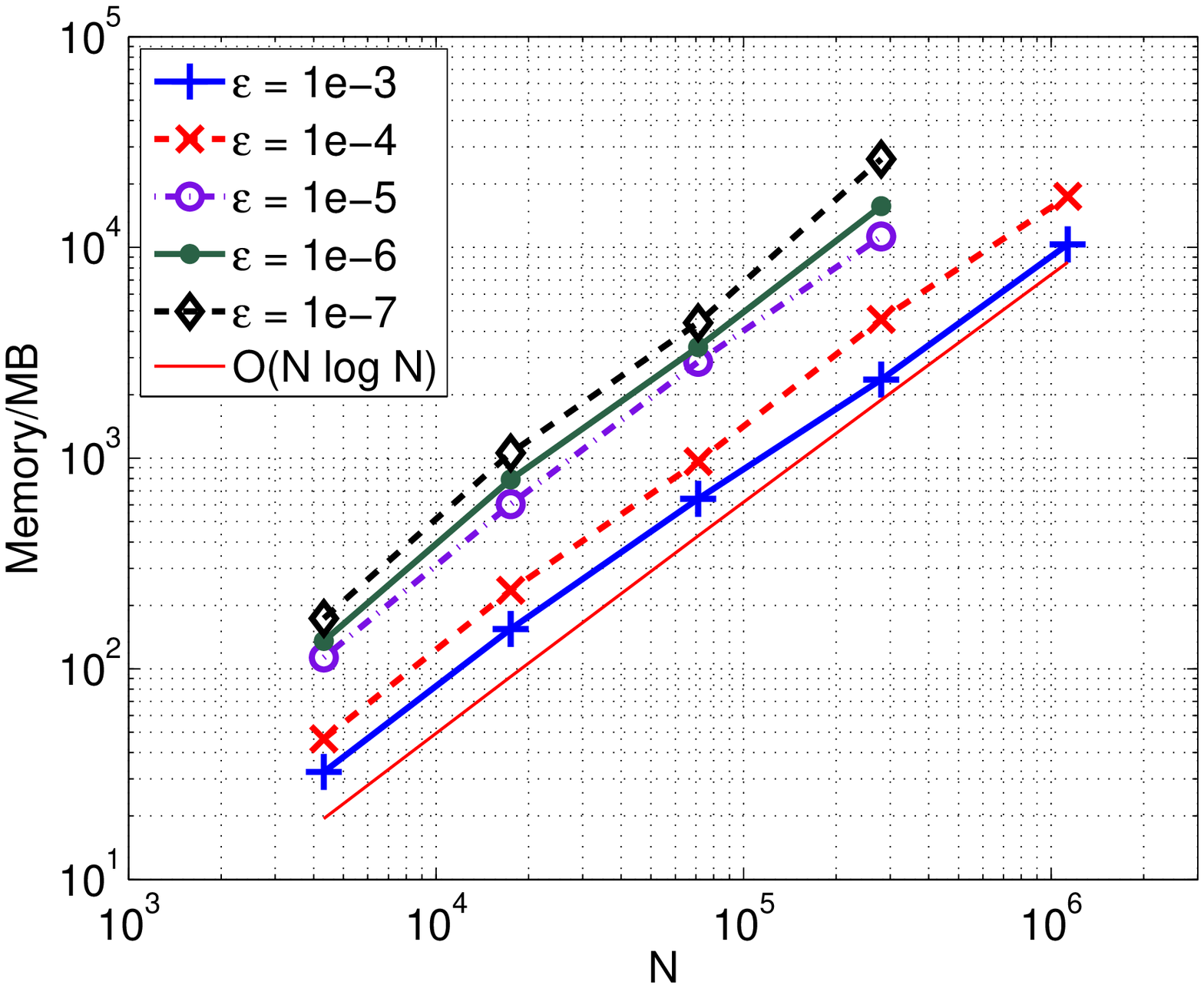}
		\label{fig_membm}
	}
	\subfigure[Computational time versus the DOF.]{
		\includegraphics[width=0.47\textwidth]{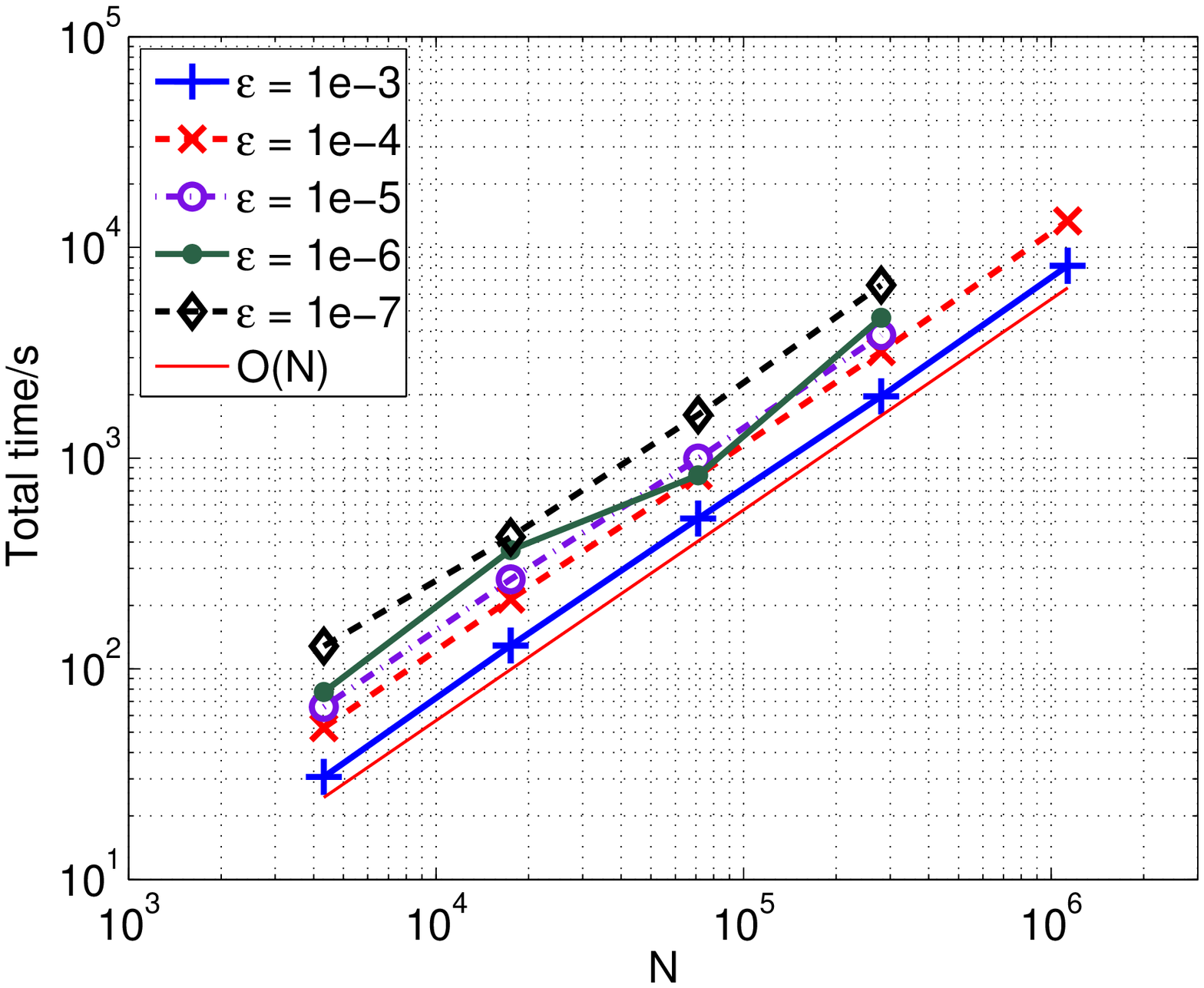}
		\label{fig_ttbm}
	}
	\subfigure[Iteration time versus the DOF.]{
		\includegraphics[width=0.47\textwidth]{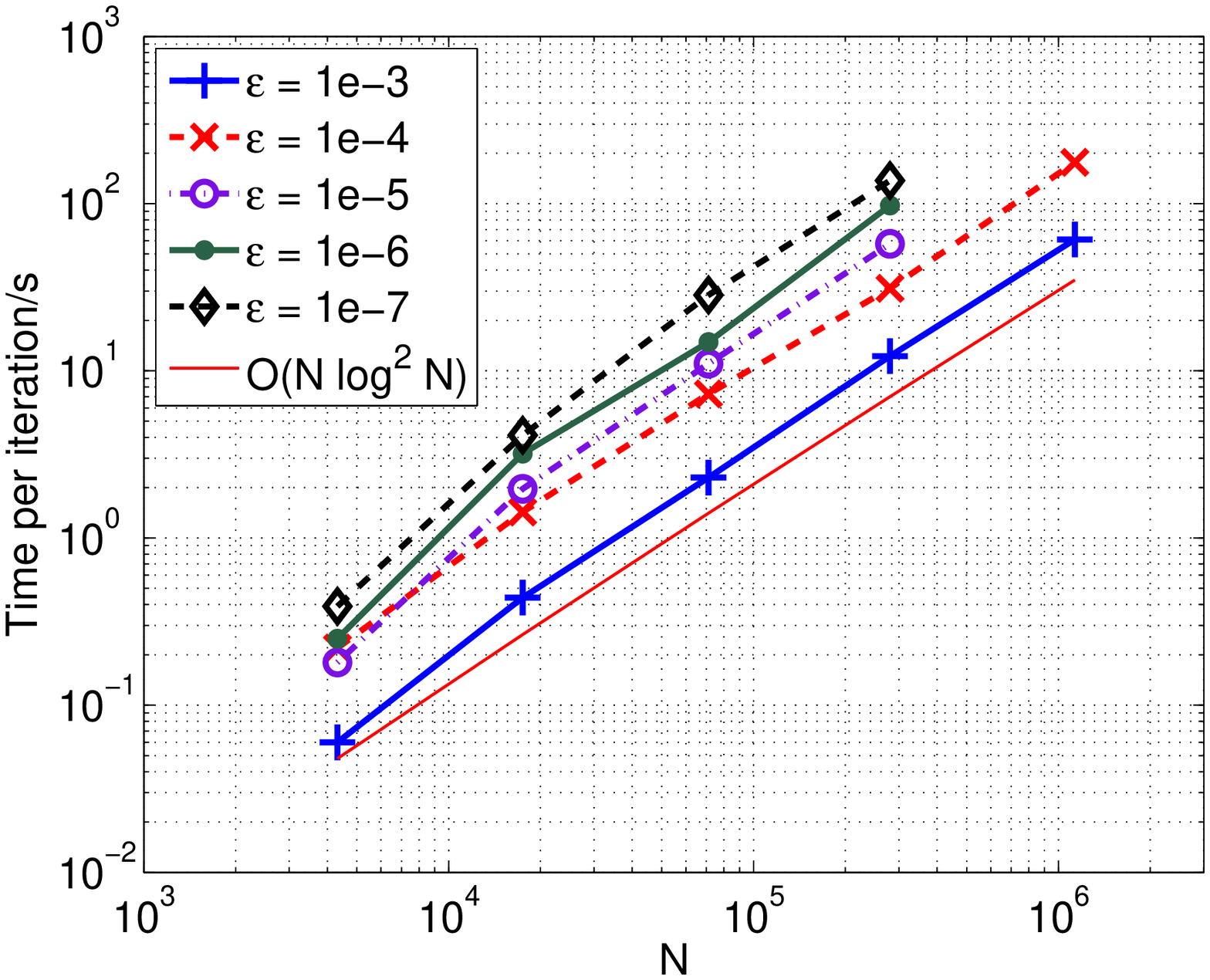}
		\label{fig_tibm}
	}
	\caption{Numerical results of the unit sphere radiating problem.}
	\label{fig_usphrad}
\end{figure}

The models are solved with different precision, determined by the
controlling accuracy $\varepsilon$, and the results are illustrated
in Figure \ref{fig_usphrad}. Figure \ref{fig_errbm} illustrates
the $L_2$ error behavior of the solution $u$ with the refinement of
mesh and the change of the controlling accuracy $\varepsilon$. It is
seen that the FDBEM can attain the given precision $O(\varepsilon)$
for $\varepsilon \gtrsim 10^{-6}$. When $\varepsilon = 10^{-7}$, the
error of the solution keeps almost the same with that of
$\varepsilon = 10^{-6}$, because in this case the error in the
numerical evaluation of the singular and near-singular boundary
element integrations becomes dominate. Meanwhile, the error curves
first go down to $O(\varepsilon)$ with the refinement of the meshes, and
then tend to rise slightly with the further refinement of the mesh.
This is because, in the first stage the BEM discretization error is
dominate, and in the second stage the error induced by the
translations in the FDA which is of order $O(\varepsilon \log N)$
\cite{kifmm} becomes dominate.

The memory consumption, overall CPU time and the CPU time for each
iteration are illustrated in Figure \ref{fig_membm}, \ref{fig_ttbm}
and \ref{fig_tibm}, respectively. 
For each precision $\varepsilon$, the CPU time for each iteration
and the total memory usage of the FDBEM are of order $O(N \log N)$.
The overall CPU time in \ref{fig_ttbm} grows almost linearly with
$N$. This is because a large part of the total CPU time is spent
on the numerical evaluation of the various singular integrals
in the near-field matrix which scales linearly with $N$. For
example, in the case with $\varepsilon =$1e-3, $k=32\pi$ and
$N=1132776$, the overall CPU time is 8181s, while evaluating
the near field part of the matrices $\bm{H}$ and $\bm{G}$
takes 7053s, which is about 86\% of the overall CPU time.

\subsection{Large-scale simulations}

\begin{figure}[h]
	\centering
	\subfigure[The real part of total acoustic field.]{
		\includegraphics[height=0.4\textwidth]{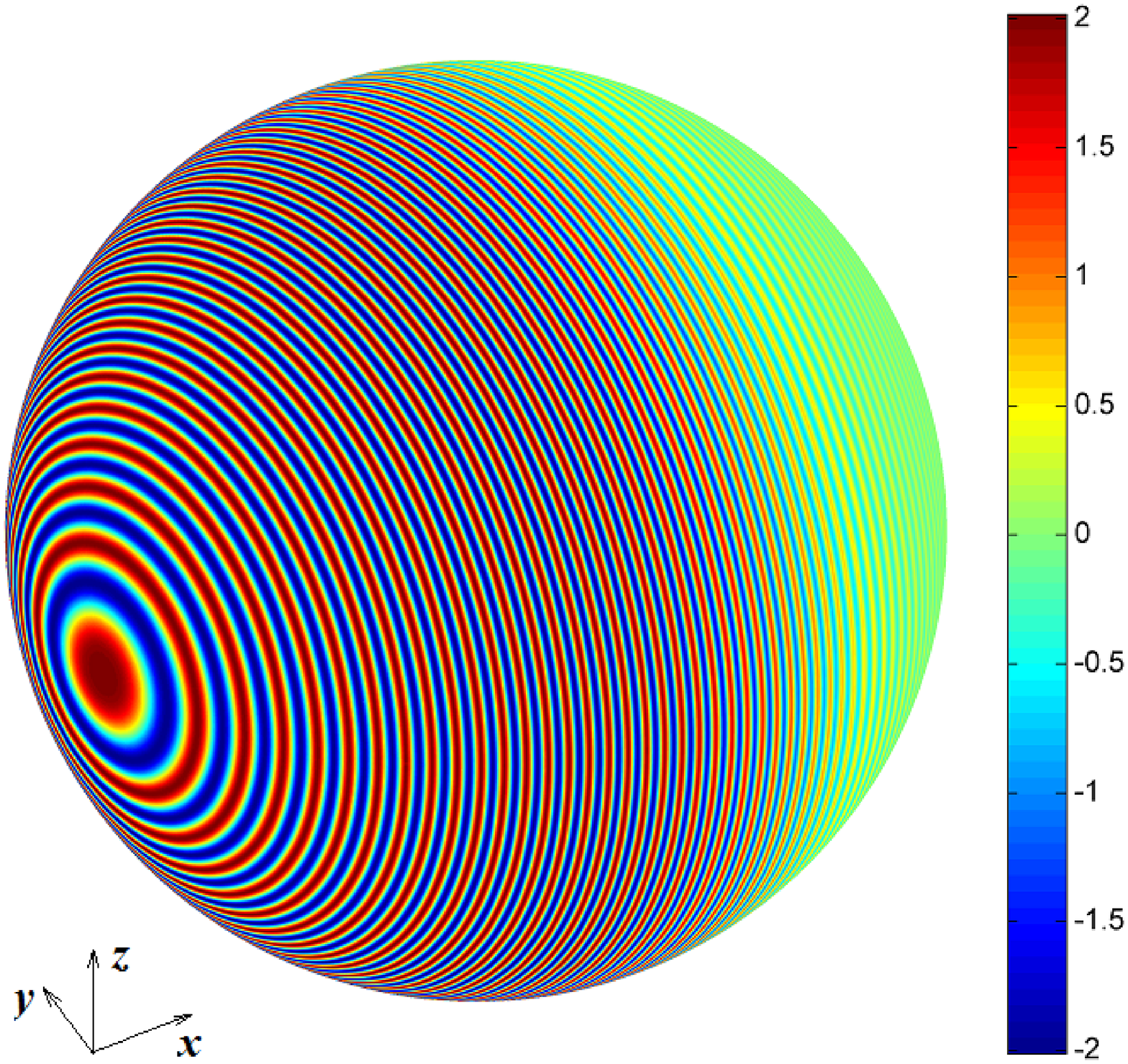}
		\label{fig_usphfield}
	}
	\subfigure[The amplitude of the scattering field.]{
		\includegraphics[height=0.4\textwidth]{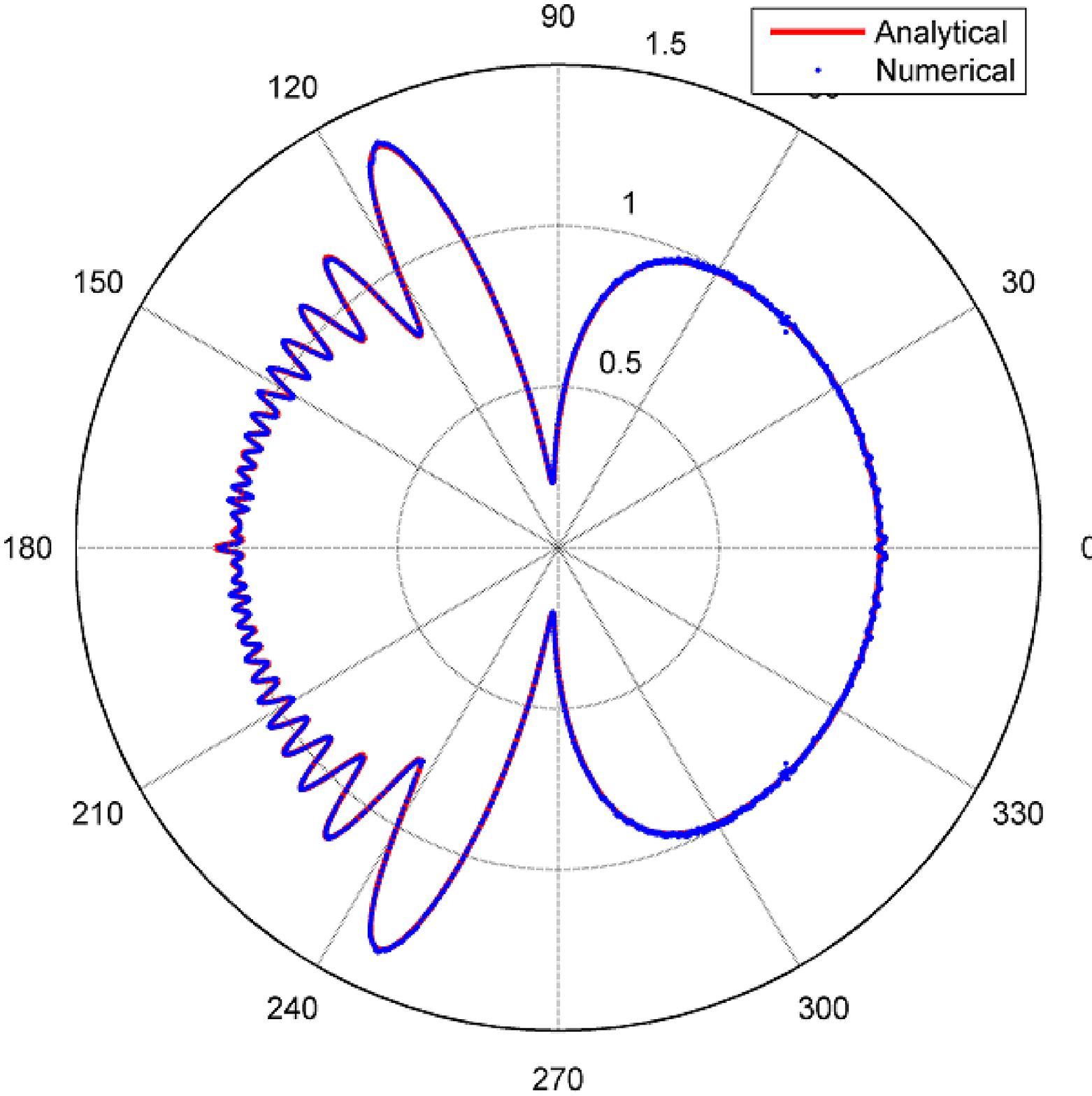}
		\label{fig_usphscaptn}
	}
	\caption{Results of the unit sphere scattering problem.}
	\label{fig_usphsca}
\end{figure}

\begin{figure}[h]
	\centering
	\includegraphics[width=0.7\textwidth]{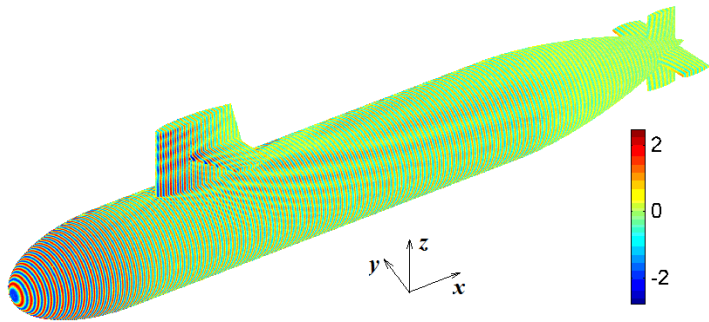}
	\caption{The real part of the total acoustic field of the submarine
		scattering problem.}
	\label{fig_submarine}
\end{figure}

Finally, the performance of the FDBEM for large-scale simulations is
demonstrated by solving two sound scattering problems. In both 
problems,
the scatters are assumed to be sound hard, so the boundary condition is 
$q = \frac{\partial u}{\partial \bm{n}} = 0$. The incident wave is a 
plane wave propagating in $\bm{d} = (1, 0, 0)$ direction, with the 
acoustic velocity potential being $u^\text{inc} (\bm{x}) = e^{\text{i}k
bm{x} \cdot \bm{d}}$, see Figure \ref{fig_usphfield} and 
\ref{fig_submarine} for the coordinate system.

The scatter for the first problem is a unit sphere with its diameter $D$ 
of 50 wavelengths long, that is, the wave number is $k = 157.1$, and the 
non-dimensional wave number is $kD = 314.2$. The analytical solution of 
the scattered acoustic velocity potential on the spherical surface 
is given by
\begin{equation}
	u^\text{s} (\theta) = \sum_{m=0}^\infty \left[ -\text{i}^m (2m+1) 
	\frac{j'_m(k)}{h'_m(k)} \right] h_m(k) \cdot P_m(\cos \theta),
\end{equation}
where $\theta$ is the angle between $\overrightarrow{O\bm{x}}$ and 
$\bm{d}$, $j_m$ is the spherical Bessel function of the first kind, 
$h_m$ is the spherical Hankel function of the second kind, and 
$P_m$ is the Legendre polynomial of order $m$. Then the analytical 
solution can be achieved by $u = u^\text{inc} + u^\text{s}$.

The spherical surface is 
discretized into 668352 curved triangular quadratic 
elements, and the DOFs $N=4040112$. The controlling accuracy and the 
GMRES converging tolerance are set to be $\varepsilon = 10^{-3}$.
It takes about 8.1 hours to solve this problem. The CPU time for each 
iteration is about 243.7 seconds. The memory usage is 27.7 GB. The 
GMRES solver converged 
within $N_\text{it} = 14$ iterations. The acoustic velocity potential
of the total sound field $u$ on the surface is illustrated in 
Fig. \ref{fig_usphfield}, and the comparison between the analytical 
solution and our numerical results of the scattered field $u^\text{s}$
on the surface is illustrated in Fig. \ref{fig_usphscaptn}. It is shown
that our numerical results agree very well with the analytical solution.

The scatter in the second problem is a submaine model with length $D$=80 
meters. The wave number is set to be $k=12.5$, thus the non-dimensional 
wave number $kD=1000$ and the length of the submarine is 159 wavelengths.
The surface of the submarine is discretized into 670564 curved
triangular quadratic elements, so the DOFs $N=4023384$. The controlling
accuracy and the GMRES converging tolerance are set to be $\varepsilon=
10^{-3}$. It takes 9.3 hours to solve this 
problem. The CPU time for each iteration is 250 seconds. The GMRES solver
converges within 42 iterations. The total memory usage is 24 GB. The 
acoustic velocity potential of the total field $u$ on the surface is 
illustrated in Fig. \ref{fig_submarine}.

\section{Conclusion}

The fast directional algorithm (FDA) recently developed in \cite{fda} is 
a highly efficient method for the evaluation of potential summations with 
Helmholtz kernels. In this paper, a FDA accelerated BEM, denoted by FDBEM, 
for solving Burton-Miller BIE in a broad frequency range is developed and 
implemented. The Nystr\"{o}m method based on curved quadratic elements is 
used to discretize the BIE, by which: (1) the edge and corner problems 
are completely avoided, (2) high accuracy can be achieved and (3) the 
resulting linear system is more like a summation that is suitable for FDA 
acceleration. The FDA for the potential summations with the kernels of the 
Burton-Miller formulation are proposed. The computational efficiency
of the FDA is further elevated by exploiting the low-rank property
of the translation matrices.

By using the FDBEM, large-scale wideband acoustic problems can be solved
with controllable accuracy up to $10^{-6}$. The computational time and
memory requirements are of order $O(N\log N)$. An representative acoustic
scattering problem with dimensionless wave number $kD$ being up to 1000 
and the DOF being up to 4 million has been successfully solved within 
10 hours on a computer with one core and the memory usage is 24 GB.

\section*{Acknowledgements}

This work is supported by the National Science
Foundation of China under Grants 11074201 and 11102154, the Funds
for Doctor Research Programs from the Chinese Ministry of
Education under Grants 20106102120009 and 20116102110006, and the 
Doctorate Foundation of Northwestern Polytechnical University under Grant
No. CX201220.

\bibliographystyle{unsrt}
\bibliography{fda}

\end{document}